\documentclass[technote,onecolumn]{IEEEtran}
\usepackage{graphicx}
\usepackage{amsfonts}
\usepackage[cmex10]{amsmath}
\usepackage[ruled]{algorithm}

\usepackage{url}

\begin{document}


\title{\begin{large}A Unifying Approach to Quaternion Adaptive Filtering: Addressing the Gradient and Convergence\end{large}}

\author{
\textbf{Cyrus Jahanchahi} and \textbf{Danilo P. Mandic}\\
Communications and Signal Processing Research Group,\\Department of Electrical and Electronic Engineering,\\ Imperial College London, London, SW7 2AZ, U.K.\\
email: \{cyrus.jahanchahi, d.mandic\}@imperial.ac.uk\\
URL: \url{http://www.commsp.ee.ic.ac.uk/~mandic}\\
\vspace{2mm}
\begin{large}
\textbf{Technical Report}: TR-IC-Quaternion-SP-01-08-2012\\
\textit{Created}: January 2012, \textit{Updated}: September 2012
\end{large}
\thanks{\textbf{\copyright~C. Jahanchahi and D. P. Mandic} --- This research is supported by EPSRC under the grant EP/H026266/1. {\bf - Personal use of this material is permitted. Permission to use this material for any other purposes must be obtained.}}
}

\maketitle

\begin{abstract}
A novel framework for a unifying  treatment of quaternion valued adaptive filtering algorithms is introduced. This is achieved based on  a rigorous account of quaternion differentiability, the proposed I-gradient, and the use of augmented quaternion statistics to account for real world data with noncircular probability distributions. We first provide an elegant solution for the calculation of the gradient of real functions of quaternion variables (typical cost function), an issue that has so far prevented systematic development of quaternion adaptive filters. This makes it possible to unify the  class of existing and proposed quaternion least mean square (QLMS) algorithms, and to illuminate their structural similarity. Next, in order to cater for both circular and noncircular data, the class of widely linear QLMS (WL-QLMS) algorithms is introduced and the subsequent convergence analysis unifies  the treatment of strictly linear and widely linear filters, for both proper and improper sources. It is also shown that the proposed class of $\mathbb{HR}$ gradients allows us to resolve the uncertainty owing to the noncommutativity of quaternion products, while the involution gradient (I-gradient) provides generic extensions of the corresponding real- and complex-valued adaptive algorithms, at a reduced computational cost. Simulations in both the strictly linear and widely linear setting support the approach.
\end{abstract}


\section{Introduction}

Recent developments in sensor technology, human centered computing, and robotics \cite{Biamino_RobotEye_2005}\cite{Miron_QuaternionMUSIC_2006} \cite{Sabatini} have made possible the recording og  new classes of multidimensional signals which are naturally represented as two- and three-dimensional vector-valued processes. Such signals are readily modeled as real vectors (in $\mathbb{R}^2$ and $\mathbb{R}^3$), however, there are big advantages in processing such data in division algebras - the complex and quaternion domain. The benefits of complex domain modeling are well understood \cite{Mandic_Book_ComplexFilters}, while  for three- and four-dimensional data quaternions (a skew field over $\mathbb{R}^4$), form a noncommutative division algebra  and provide a unifying treatment of two-, three- and four-dimensional processes.

Additional benefits arising from the use of such a 3D and 4D division algebra include very accurate rotation and orientation modeling. This has made quaternions a standard in computer graphics\cite{Shoemake_SLERP}, where they  have replaced real rotation matrices, owing to their ability to  provide a compact and accurate representation of sequences of rotations, together with  a natural account of the spherical angle \cite{Hanson_Visulazing_Quaternions}. This has brought new opportunities in a variety of applications, such as color imaging \cite{Pei_ColorImageProcessing_1999}, body-sensor networks \cite{Rehman_Empirical_Mode_Decomposition_2009}, and aeronautics \cite{Choukroun_QuaternionKalman_2006}.

Although the ordering relationships (e.g. `$\geq$') are not defined for quaternions (and hence the notion of probability density function), there exists a `static' isomorphism between $\mathbb{R}^3$ and $\mathbb{R}^4$ and the field of quaternions $\mathbb{H}$ (pure and full), which can be used  to introduce  the notions of probability distribution and gradient. Such tools are a prerequisite for the design of learning algorithms in the quaternion domain, which offers theoretical and practical advantages when dealing with the large (and rotational) dynamics of vector signals. Examples  are quite recent and include the quaternion least mean square (QLMS) \cite{Clive_QLMS_2009} and  developments in `augmented' quaternion statistics \cite{Clive_Statistics_2010}  \cite{Bihan_quaternion_properness} \cite{Via_properness}. These  have been followed by \emph{widely linear} adaptive filtering algorithms, including the widely linear QLMS (WL-QLMS) \cite{Clive_WL-QLMS_2009} and widely linear blind source separation algorithms for improper quaternion sources \cite{Soroush_BSS} \cite{Soroush_FastICA} \cite{Via_ICA}.

Standard, or strictly linear, statistical models rest upon the covariance matrix $E[\mathbf{q}\mathbf{q}^H]$, where $\mathbf{q} \in \mathbb{H}^{N \times 1}$ is a regressor vector, and are  second order optimal only for circular signals, for which the probability distributions are rotation invariant. In terms of second order statistics, such signals have equal powers in all the four quaternion components and are termed `proper', however, real world processes are typically `improper' (with different component powers). To fully describe their second order properties, three complementary covariance matrices must also be employed  \cite{Clive_Statistics_2010} \cite{Via_properness}.

 A rigorous treatment of the gradient is a key prerequisite to the derivation of stochastic gradient optimization procedures, and the inconsistencies and lack of coherence in their definition has so far prevented a more widespread use of quaternions. Current approaches typically employ quaternion gradient calculation based on the so called `pseudogradient' \cite{Wieslaw_Fueter_2006}, however, this gradient, also known as the Cauchy-Riemann-Fueter (CRF) derivative, is very stringent, and poses a number of restrictions in the derivation of the adaptive filtering algorithms. For instance, due to the ambiguity related to the noncommutativity of quaternions, the CRF derivative can be calculated with the imaginary unit vectors $\imath$, $\jmath$ and $\kappa$ on either the left or  right hand side of the gradient components, introducing uncertainty in choosing a correct form. Another problem with the CRF derivative is its uniqueness, as exemplified on the linear function $J=w$, whereby we obtain $\frac{\partial J}{\partial w^*}=0$ when  differentiating directly with respect to  $w^*$ and $\frac{\partial J}{\partial w^*} = -\frac{1}{2}$ using the CRF derivative. These inconsistencies have led to numerous different expressions for structurally similar algorithms, and have slowed down the development and applications of quaternion valued statistical signal processing.

The recent progress in the statistics of quaternion variable (augmented quaternion statistics \cite{Via_properness} \cite{Clive_Statistics_2010}) has facilitated the development of advanced algorithms for vector sensor modeling, however, this has also highlighted   several key remaining issues that need to be addressed in order to exploit the full power of quaternions in statistical signal processing applications. To that end, this paper address the following problems:

\begin{itemize}
  \item Unifying different forms of QLMS. Due to the noncommutativity of the quaternion product numerous expressions for the quaternion LMS exist. Our analysis  shows that the different forms are equivalent in the way they use the available information.
  \item Generic extensions of LMS. The existing quaternion LMS algorithms are not a generic extension of the real and complex LMS (CLMS). By introducing a new gradient operator (the I-gradient) we provide a generic expression for the quaternion LMS  (the IQLMS).
  \item Unifying convergence analysis. An in depth analysis for the convergence in the mean of the strictly linear QLMS is still lacking. To that end, we provide a convergence analysis for the two  existing  variants of the QLMS and for the novel IQLMS. The analysis  shows that the IQLMS offers the same performance  while requiring only half the mathematical operations.
  \item A rigorous MSE analysis. No analysis exists for the convergence in the mean square of the widely linear QLMS algorithms. We provide a rigorous and unifying mean square convergence analysis  which quantifies the advantage offered by the widely linear algorithms  for noncircular signals.
\end{itemize}
In this way, we have provided a unifying set of tools and an analysis platform for future developments in the field, which has many theoretical benefits but has found only a few applications due to the above uncertainties.
The generic form of algorithms enabled by the I-gradient makes it possible to extend the vast resources from the real and complex adaptive filtering into the field of quaternions.

 The organization of the paper is as follows: in Section \ref{sec:quat_algebra}, we briefly review the elements of quaternion algebra necessary for the development of the $\mathbb{HR}$-calculus and quaternion valued filters presented. In Section \ref{Sec:HRcalculus}, the $\mathbb{HR}$-calculus is summarized   and it is shown how this theory unifies  gradient calculations. Section \ref{Sec:WidelyLinearModel} introduces quaternion valued widely linear models. In Section \ref{Sec:QuaternionGradient}, the conjugate gradient is addressed and a low complexity gradient expression (the I-gradient) is proposed.  In Section \ref{Sec:Unified_approach}, the $\mathbb{HR}$-calculus is  used to efficiently re-derive the QLMS and WL-QLMS algorithms in a unifying manner and to introduce a novel class of algorithms based on the I-gradient, requiring half the mathematical operations of the QLMS while having the same generic form as LMS and CLMS. In Sections \ref{Section:equivalence} and \ref{Sec:Convergence}, we analytically compare the performances of the IQLMS and QLMS algorithms, introducing a generic framework for the analysis of quaternion valued adaptive filters. In Section \ref{sec:simulations}, the performances of all algorithms considered are illustrated for both circular (proper) and second order noncircular (improper signals).


\section{Quaternion Algebra \label{sec:quat_algebra}}
Quaternions are a skew field over $\mathbb{R}^4$ defined as
\vspace{0mm}
\begin{equation*}
\{q_r,q_{\imath},q_{\jmath},q_{\kappa}\} \in \mathbb{R}^4 \rightarrow q= q_r + {\imath}q_{\imath}  + {\jmath}q_{\jmath}  + {\kappa}q_{\kappa} \in \mathbb{H}
\vspace{0mm}
\end{equation*}
 Unlike in $\mathbb{R}^4$ the unit axis vectors $\imath$, $\jmath$ and $\kappa$ are also the imaginary units,  and obey the following rules
 \vspace{0mm}
\begin{equation*}
\imath \jmath = \kappa \text{ } \text{ } \text{ } \jmath \kappa = \imath  \text{ } \text{ } \text{ } \kappa \imath = \jmath \nonumber
\end{equation*}
\begin{equation*}
\imath^2 = \jmath^2 = \kappa^2 = \imath\jmath\kappa = -1
\vspace{0mm}
\end{equation*}
Note that quaternion multiplication is not commutative, that is, $\imath \jmath \neq \jmath\imath = -\kappa$.

A quaternion variable $q \in \mathbb{H}$ can be more conveniently expressed as \cite{Ward_Book_QuaternionsCayleyNumbers}
\vspace{0mm}
\begin{equation*}
q= Sq + Vq
\vspace{0mm}
\end{equation*}
where $Sq= q_r $ denotes the scalar part of $q$ and $Vq= \imath q_{\imath}  + {\jmath}q_{\jmath}  + {\kappa}q_{\kappa}$ the vector part. The product of quaternions $q_1,q_2 \in \mathbb{H}$ is given by
\vspace{0mm}
\begin{eqnarray*}
q_1q_2 &=& (Sq_1 + Vq_1)(Sq_2 + Vq_2)
= Sq_1Sq_2 - Vq_1 \bullet Vq_2
+ Sq_2Vq_1
 + Sq_1Vq_2 + Vq_1 \times Vq_2
\end{eqnarray*}
where the symbol `$\bullet$' denotes the dot-product and `$\times$' the cross-product.
The quaternion conjugate, denoted by $q^*$, is defined  as
\vspace{0mm}
\begin{equation}
q^*= Sq -Vq = q_r - \imath q_{\imath} - \jmath q_{\jmath} - \kappa q_{\kappa} \vspace{0mm}
\end{equation}
and the norm $\parallel q \parallel$   as
\vspace{0mm}
\begin{equation*}
\parallel q \parallel= \sqrt{qq^*}=\sqrt{q^2_{r} + q^2_{\imath} + q^2_{\jmath} +q^2_{\kappa}}
\vspace{0mm}
\end{equation*}
The three-dimensional vector part  $Vq$ is also called a \emph{pure quaternion}, while the inclusion of the real part $Sq$ gives a \emph{full quaternion}. Notice that the same notation $Vq$ is used for describing both a pure quaternion and a three dimensional real vector. The difference is in that when describing a vector in $\mathbb{R}^3$, $\{\imath,\jmath,\kappa\}$ are  unit vectors, whereas for a pure quaternion, $\{\imath,\jmath,\kappa\}$ are not only unit vectors but also imaginary units. In this way, the algebraic structure of quaternions enables a unified processing of both three- and four-dimensional vector processes.

\subsection{Quaternion Involutions}
Quaternion involutions  are similarity relations, or self-inverse mappings, defined as\footnote{Note that the quaternion conjugate is also an involution, that is, a self-inverse mapping.} \cite{involutions}
\vspace{0mm}
\begin{eqnarray}
q^{\imath} &=& -\imath q \imath = q_r + {\imath}q_{\imath} - \jmath q_{\jmath} - \kappa q_{\kappa}  \nonumber \\
q^{\jmath} &=& -\jmath q \jmath = q_r - \imath q_{\imath}  + {\jmath}q_{\jmath} - \kappa q_{\kappa} \nonumber \\
q^{\kappa} &=& -\kappa q \kappa  = q_r - \imath q_{\imath} - \jmath q_{\jmath}  + {\kappa}q_{\kappa}
\label{eq:involution}
\vspace{0mm}
\end{eqnarray}
 and have the following properties (for $inv_3 \neq inv_2 \neq inv_1$)
\begin{eqnarray}
   P1: &&( q ^{inv} )^{inv} = q \text{ for } inv \in \{\imath ,\jmath ,\kappa \} \hspace{15mm}
   P2: ( q_1 q_2 )^{inv} = q_1^{inv} q_2^{inv} \nonumber \\
   P3: &&( q_1 + q_2 )^{inv} = q_1^{inv} + q_2^{inv} \hspace{24mm}
   P4:  (q^{inv_1})^{inv_2} =  (q^{inv_2})^{inv_1} = q^{inv_3}
   \label{eq:algebra:inv:p4}
\end{eqnarray}
 Involutions can be seen as a generalisation of the complex conjugate,  as they allow for the components of a quaternion variable $q$ to be expressed in terms of the actual variable $q$ and its `partial conjugates' $q^{\imath}$, $q^{\jmath}$, $q^{\kappa}$, that is\footnote{Compare with the complex domain where the real and imaginary parts of the complex numbers $z=x+{\imath}y$ are expressed by
$
x= \frac{1}{2}(z + z^*) \text{ and } y= \frac{1}{2i}(z - z^*)
$.}
\begin{eqnarray}
q_r &=&  \frac{1}{4} [q + q^{\imath} + q^{\jmath} + q^{\kappa}] \hspace{10mm}  q_{\imath} =  \frac{1}{4 \imath} [q + q^{\imath} - q^{\jmath} - q^{\kappa}] \nonumber \\
q_{\jmath} &=&  \frac{1}{4 \jmath} [q - q^{\imath} + q^{\jmath} - q^{\kappa}] \hspace{10mm}  q_{\kappa} =  \frac{1}{4 \kappa} [q - q^{\imath} - q^{\jmath} + q^{\kappa}]
\label{Eq:QuatAlgebra:q_d}
\end{eqnarray}

This representation will be instrumental in the derivation of both the $\mathbb{HR}$-calculus and the quaternion \emph{widely linear} model, and underpins the development of  quaternion valued adaptive filters for noncircular signals.
\subsection{Quaternion Algebra Versus Three-Dimensional Vector Algebra}

Consider a transformation matrix $\mathbf{A} \in \mathbb{R}^{3 \times 3}$ that maps a point $\mathbf{x} \in \mathbb{R}^{3}$ onto a point $\mathbf{y} \in \mathbb{R}^{3}$ and a quaternion $q_A \in \mathbb{H}$ that transforms a pure quaternion $q_x \in \mathbb{H}$ into  a pure quaternion $q_y \in \mathbb{H}$, that is
\vspace{0mm}
\begin{equation}
\mathbf{y}=\mathbf{Ax} \hspace{5mm} \in \mathbb{R}^{3} \hspace{1cm} \sim \hspace{1cm} q_y=q_A q_x q_A^* \hspace{5mm}\in \mathbb{H}
\label{eq:quat_vs_vectors}
\vspace{0mm}
\end{equation}

\noindent \textbf{Remark\#1:}
The mapping $\mathbf{A}$ requires nine coefficients to relate two vectors in $\mathbb{R}^3$, however, physically only four parameters are needed (two for the axis of rotation, one for the angle of rotation and one for the scaling factor). The four elements of a quaternion offer  this physical insight and compact representation, expressing straightforwardly  the  axis of rotation, angle of rotation, and scaling factor.


\noindent \textbf{Remark\#2:}
When the mapping $\mathbf{A}$ represents a succession of rotations in the $x$, $y$, $z$ directions (using Euler angles), a degree of freedom can be lost if any two axis are aligned, resulting in the gimbal lock phenomenon. This is not a problem in the quaternion domain, where the quaternion transformation in (\ref{eq:quat_vs_vectors}) is expressed as $q_y=q_A q_xq^{-1}_A$, where $q_A$ is a unit quaternion.

\noindent \textbf{Remark\#3:}
The  quaternion rotation $q_A$ is better conditioned than the real rotation matrix $\mathbf{A}$, as the only requirement on $q_A$ is to be  unit quaternion, whereas $\mathbf{A}$ must satisfy $\mathbf{A}^T\mathbf{A}=\mathbf{I}$ and $det(\mathbf{A})=1$. This  has led to the use of quaternions in e.g. spacecraft orientation problems where they provide convenient closed form solutions \cite{Choukroun_QuaternionKalman_2006} \cite{Marins_ExtendedKalman_Mars} \cite{Horn_closed_formsolution}.

\section{The $\mathbb{HR}$-Calculus \label{Sec:HRcalculus}}

In adaptive filtering, the aim  of gradient based optimization is to minimise a measure of error power, typically a scalar function of quaternion variables. This has been a main stumbling block in the development of learning algorithms, as the Cauchy-Riemann-Fueter (CRF) conditions
 \begin{equation}
\frac{\partial J}{\partial \mathbf{w}^*}= \frac{1}{4} \left(\frac{\partial J}{\partial \mathbf{w}_r} + {\imath}\frac{\partial J}{\partial \mathbf{w}_{\imath}}  + {\jmath}\frac{\partial J}{\partial \mathbf{w}_{\jmath}}  + {\kappa}\frac{\partial J}{\partial \mathbf{w}_{\kappa}} \right)=\mathbf{0}
\label{Eq:HR:Fueter}
\end{equation}
 do not admit the calculation of such gradients since they are not defined for real functions. To that end, following on the corresponding result in the complex domain,  where the $\mathbb{CR}$-calculus \cite{Mandic_Book_ComplexFilters} \cite{Delgado_Report_ComplexGradient} is used to calculate the gradient of real functions of complex variables, we have recently introduced the $\mathbb{HR}$-calculus \cite{Cyrus_SPL},
which consists of two groups of derivatives: the $\mathbb{HR}$-derivatives
\begin{equation}
\left[ \begin{array}{cccc}
\frac{\partial f(q,q^\imath,q^\jmath,q^\kappa)}{\partial q} \\
\frac{\partial f(q,q^\imath,q^\jmath,q^\kappa)}{\partial q^\imath} \\
\frac{\partial f(q,q^\imath,q^\jmath,q^\kappa)}{\partial q^\jmath} \\
\frac{\partial f(q,q^\imath,q^\jmath,q^\kappa)}{\partial q^\kappa} \\
\end{array} \right]= \frac{1}{4}
\left[ \begin{array}{cccc}
 1 & -\imath & -\jmath & -\kappa\\
 1 &  -\imath & \jmath & \kappa\\
 1 &  \imath &  -\jmath & \kappa\\
 1 &  \imath &  \jmath & -\kappa \end{array} \right]
 \left[ \begin{array}{cccc}
\frac{\partial f(q_r,q_{\imath},q_{\jmath},q_{\kappa})}{\partial q_r} \\
\frac{\partial f(q_r,q_{\imath},q_{\jmath},q_{\kappa})}{\partial q_i} \\
\frac{\partial f(q_r,q_{\imath},q_{\jmath},q_{\kappa})}{\partial q_{\jmath}} \\
\frac{\partial f(q_r,q_{\imath},q_{\jmath},q_{\kappa})}{\partial q_{\kappa}} \\
\end{array} \right]
\label{Eq:HRcalculus:HRderivatives}
\end{equation}
%
%
and the $\mathbb{HR}^*$-derivatives
\begin{equation}
\left[ \begin{array}{c}
\frac{\partial f(q^*,q^{\imath\!*},q^{\jmath\!*},q^{\kappa\!*})}{\partial q^*} \\
\frac{\partial f(q^*,q^{\imath\!*},q^{\jmath\!*},q^{\kappa\!*})}{\partial q^{\imath\!*}} \\
\frac{\partial f(q^*,q^{\imath\!*},q^{\jmath\!*},q^{\kappa\!*})}{\partial q^{\jmath\!*}} \\
\frac{\partial f(q^*,q^{\imath\!*},q^{\jmath\!*},q^{\kappa\!*})}{\partial q^{\kappa\!*}} \\
\end{array} \right]= \frac{1}{4}
\left[ \begin{array}{cccc}
 1 & \imath & \jmath & \kappa\\
 1 &  \imath & -\jmath & -\kappa\\
 1 &  -\imath &  \jmath & -\kappa\\
 1 &  -\imath &  -\jmath & \kappa \end{array} \right]
 \left[ \begin{array}{cccc}
\frac{\partial f(q_r,q_{\imath},q_{\jmath},q_{\kappa})}{\partial q_r} \\
\frac{\partial f(q_r,q_{\imath},q_{\jmath},q_{\kappa})}{\partial q_i} \\
\frac{\partial f(q_r,q_{\imath},q_{\jmath},q_{\kappa})}{\partial q_{\jmath}} \\
\frac{\partial f(q_r,q_{\imath},q_{\jmath},q_{\kappa})}{\partial q_{\kappa}} \\
\end{array} \right]
\label{Eq:HRcalculus:HR*derivatives}
\end{equation}
where $q^{inv*}=(q^{inv})^* \text{ for } inv \in \{\imath ,\jmath ,\kappa \}$ are the conjugates of the involutions in (\ref{eq:involution}).

\noindent \textbf{Remark\#4:}
The $\mathbb{HR}^*$-derivative $\frac{\partial f(q^*,q ^{{\imath}*},q ^{{\jmath}*},q ^{{\kappa}*})}{\partial q^*}$ is equivalent to the quaternion derivative operator introduced by Fueter \cite{Wieslaw_Fueter_2006}, however, unlike the CRF derivative in (\ref{Eq:HR:Fueter}), the derivative $\frac{\partial f(q^*,q ^{{\imath}*},q ^{{\jmath}*},q ^{{\kappa}*})}{\partial q^*}$ also imposes a restriction on the form of dependent variables that compose the function $f(\cdot)$. This, for the first time, allows for a direct differentiation of quaternion functions, similar to that in $\mathbb{R}$ and $\mathbb{C}$.

\noindent \textbf{Remark\#5:}
The $\mathbb{HR}$-derivatives and $\mathbb{HR}^*$-derivatives can be used in a similar way to the $\mathbb{R}$- and $\mathbb{R}^*$-derivative in the complex domain \cite{Delgado_Report_ComplexGradient}. For instance, to perform a direct $\mathbb{HR}$ differentiation of a function written in terms of $q^*$, we   must first rewrite it in terms of $q$, $q^{\imath}$, $q^{\jmath}$ and $q^{\kappa}$, using the substitution
\vspace{0mm}
\begin{equation}
q^*=\frac{1}{2}(q^{\imath} + q^{\jmath} + q^{\kappa} - q)
\label{Eq:UseOfHR:q*}
\vspace{0mm}
\end{equation}
Similarly, to differentiate a function  of $q$ using the $\mathbb{HR}^*$-derivatives, we substitute for $q$ using
\vspace{0mm}
\begin{equation}
q=\frac{1}{2}(q ^{{\imath}*} + q ^{{\jmath}*} + q ^{{\kappa}*} - q^*)
\label{Eq:UseOfHR:q}
\vspace{0mm}
\end{equation}

 Therefore, the $\mathbb{HR}$-calculus provides a tool for differentiating quaternion functions directly, rather than employing partial derivatives with respect to the real valued $q_r$, $q_{\imath}$, $q_{\jmath}$, $q_{\kappa}$, as is current practice (within the pseudogradient).

\noindent \textbf{Remark\#6:}
The derivatives in (\ref{Eq:HRcalculus:HRderivatives}) and (\ref{Eq:HRcalculus:HR*derivatives}) have the imaginary unit vectors placed on the left hand side of the real partial derivatives $\frac{\partial f}{\partial q_r}$, $\frac{\partial f}{\partial q_{\imath}}$, $\frac{\partial f}{\partial q_{\jmath}}$ and $\frac{\partial f}{\partial q_{\kappa}}$ (termed the `left-$\mathbb{HR}$ and left-$\mathbb{HR}^*$-derivatives'). Due to the noncommutativity of the quaternion product, an alternative set of derivatives can be obtained by writing (\ref{Eq:HRcalculus:HRderivatives}) and (\ref{Eq:HRcalculus:HR*derivatives}) as row vectors, with the unit vectors placed on the right hand side of the real partial derivatives, giving rise to the `right-$\mathbb{HR}$ and right-$\mathbb{HR}^*$-derivatives'. Both are equally valid \cite{Cyrus_SPL}, though only one arrangement will give the correct derivative (depending on the function differentiated). The $\mathbb{HR}$-calculus rectifies this ambiguity by directly differentiating a function with respect to a quaternion variable (see Remark 5), thus removing the uncertainty related to  the position of the imaginary unit vectors, which arises when performing component-wise differentiation.

\section{Widely Linear Quaternion Modeling \label{Sec:WidelyLinearModel}}

The existing  (strictly linear) estimation model in the quaternion domain  is given by
\vspace{0mm}
\begin{equation}
\hat{y}=\mathbf{w}^T\mathbf{x}
\vspace{0mm}
\end{equation}
Observe that for all the quaternion components\footnote{Throughout this paper, a vector  $\mathbf{x}$ and its involutions are treated formally as independent variables. This is a usual formalism inherited from the complex domain, and in the $\mathbb{CR}$-calculus \cite{Brandwood} \cite{VanDenBos_ComplexGradient_1994}.}
\vspace{0mm}
\begin{equation*}
\ \hat{y}_{\eta}=E[y_{\eta}\!\mid\! \mathbf{x}_{r},\mathbf{x}_{\imath},\mathbf{x}_{\jmath},\mathbf{x}_{\kappa}] \hspace{5mm}\eta \in \{r,{\imath},{\jmath},{\kappa}\}
\vspace{0mm}
\end{equation*}
and using the involutions in (\ref{eq:involution}), we can express the components of a quaternion via its involutions e.g. $x_r=\frac{1}{4}(x+x^{\imath} +x^{\jmath}+x^{\kappa})$, leading to
\vspace{0mm}
\begin{equation*}
\hat{y_{\eta}}=E[y_{\eta}\!\mid\! \mathbf{x},\mathbf{x}^{\imath},\mathbf{x}^{\jmath},\mathbf{x}^{\kappa}] \text{ and }
\hat{y}=E[y\!\mid\! \mathbf{x},\mathbf{x}^{\imath},\mathbf{x}^{\jmath},\mathbf{x}^{\kappa}]
\vspace{0mm}
\end{equation*}
In other words, since every quaternion component is a function of its involutions, to capture the full second order information available we can introduce the  \emph{widely linear} model
\vspace{0mm}
\begin{equation}
\hat{y}= \mathbf{u}^T\mathbf{x} + \mathbf{v}^T\mathbf{x}^{\imath} + \mathbf{g}^T\mathbf{x}^{\jmath} + \mathbf{h}^T\mathbf{x}^{\kappa} = \mathbf{w}^{aT}\mathbf{x}^a
\label{eq:WL-model}
\vspace{0mm}
\end{equation}
where the augmented coefficient vector  $\mathbf{w}^a=[\mathbf{u}^T,\mathbf{v}^T,\mathbf{g}^T,\mathbf{h}^T]^T$ and the augmented regressor vector  $\mathbf{x}^a=[\mathbf{x}^T,\mathbf{x}^{{\imath}T},\mathbf{x}^{{\jmath}T},\mathbf{x}^{{\kappa}T}]^T$. It should be mentioned that an alternative widely linear model could be obtained using a conjugate  augmented regressor vector $\mathbf{x}^a=[\mathbf{x}^H,\mathbf{x}^{{\imath}H},\mathbf{x}^{{\jmath}H},\mathbf{x}^{{\kappa}H}]^T$, for more detail see \cite{Clive_Statistics_2010}.

Current statistical signal processing in $\mathbb{H}$ is largely based on strictly linear models,  drawing upon the covariance matrix $\mathbf{R}_x= E[\mathbf{xx}^H]$. However,  to model both the second order circular (proper) and second order noncircular (improper) signals, based on the widely linear model in (\ref{eq:WL-model}) we need to employ the augmented covariance matrix, given by  \cite{Clive_Statistics_2010}
\begin{equation}
 \mathbf{R}_x^a = E[\mathbf{x}^a\mathbf{x}^{aH}]=
 \left[ \begin{array}{cccc}
 \mathbf{R}_x &\mathbf{ P}_x & \mathbf{S}_x & \mathbf{T}_x\\
 \mathbf{P}_x^{\imath} & \mathbf{R}_x^{\imath} & \mathbf{T}_x^{\imath} & \mathbf{S}_x^{\imath}\\
 \mathbf{S}_x^{\jmath} & \mathbf{T}_x^{\jmath} & \mathbf{R}_x^{\jmath} & \mathbf{P}_x^{\jmath}\\
 \mathbf{T}_x^{\kappa} & \mathbf{S}_x^{\kappa} & \mathbf{P}_x^{\kappa} & \mathbf{R}_x^{\kappa}\\
 \end{array} \right]
 \label{eq:widely_linear_Rqq}
\end{equation}
where  $\mathbf{R}_x =  E[\mathbf{x}\mathbf{x}^H]$, $\mathbf{P}_x =  E[\mathbf{x}\mathbf{x} ^{{\imath}H}]$, $\mathbf{S}_x =  E[\mathbf{x}\mathbf{x} ^{{\jmath}H}]$ and $\mathbf{T}_x =  E[\mathbf{x}\mathbf{x} ^{{\kappa}H}]$.

For proper signals, all the pseudocovariance matrices $\mathbf{P}_x$, $\mathbf{S}_x$ and $\mathbf{T}_x$ vanish, and such signals  have probability distributions that are rotation invariant with respect to all the six possible pairs of axes (combinations of $\imath$, $\jmath$ and $\kappa $) \cite{Clive_Statistics_2010} \cite{Via_properness}, and thus equal powers in all the components.

\noindent \textbf{Remark\#7:} The processing in $\mathbb{R}^4$ requires ten  covariance matrices, as opposed to four in the quaternion domain (since only $\mathbf{R}_x$, $\mathbf{P}_x$, $\mathbf{S}_x$ and $\mathbf{T}_x$ are needed to fully describe $\mathbf{R}^a_{x}$). Although using the quaternion covariance matrices  for four-channel real data does not offer performance advantages, this offers a much  more convenient representation for  the modeling of signal noncircularity.

\section{Unifying the Quaternion Gradient \label{Sec:QuaternionGradient}}
    Consider the first-order multivariate Taylor series expansion (TSE) of a function $f(q,q^{\imath},q^{\jmath},q^{\kappa}) : \mathbb{H}^4 \rightarrow \mathbb{R}$, given by
     \begin{equation}
     df = \frac{\partial f}{\partial q}dq + \frac{\partial f}{\partial q^{\imath}}dq^{\imath} + \frac{\partial f}{\partial q^{\jmath}}dq^{\jmath} + \frac{\partial f}{\partial q^{\kappa}}dq^{\kappa}
     \label{Eq:QuatGrad:df}
     \end{equation}
     For a real function of quaternion variables (e.g. the typical cost function $f(q)=ee^*=|e|^2$), we have
     \begin{eqnarray}
     \frac{\partial f^{\eta}(q)}{\partial q}&=& \left( \frac{\partial f(q)}{\partial q^{\eta}}  \right)^{\eta}
     \label{Eq:QuaternionGradient:Identity}
     \end{eqnarray}
      for $\eta = \{\imath,\jmath,\kappa \}$.
      Using the involution property P2 in (\ref{eq:algebra:inv:p4}), it follows that
        \begin{equation}
     df = \frac{\partial f}{\partial q}dq + \left(\frac{\partial f}{\partial q^{\imath}}dq \right)^{\imath} + \left(\frac{\partial f}{\partial q^{\jmath}}dq\right)^{\jmath} + \left(\frac{\partial f}{\partial q^{\kappa}}dq\right)^{\kappa}
     \label{Eq:QuatGrad:df2}
     \end{equation}
    Upon applying the identity (\ref{Eq:UseOfHR:q*})  and the conjugate property $(q_1q_2)^*=q^*_2q^*_1$, this yields
     \begin{equation*}
     df = \frac{\partial f}{\partial q}dq + 2dq^*\frac{\partial f}{\partial q^*} + \frac{\partial f}{\partial q}dq = 2\left(  \frac{\partial f}{\partial q}dq + dq^*\frac{\partial f}{\partial q^*} \right)
     \end{equation*}
     finally  arriving at\footnote{Using the relationship between a quaternion and its conjugate, that is, $q_r= \frac{1}{2}(q + q^*)$.}
     \begin{equation*}
     df = 4\Re\left( \frac{\partial f}{\partial q}dq \right)
     \end{equation*}
     where the symbol $\Re(\cdot)$ denotes the real (scalar) part of a quaternion.
     This result can be extended to the quaternion vector variables, where $\mathbf{q}= [q_1,q_2,...,q_N]^T$ and $\nabla_q f=\frac{\partial f}{\partial \mathbf{q}} = \left[ \frac{\partial f}{\partial q_1}, \frac{\partial f}{\partial q_2},..., \frac{\partial f}{\partial q_N} \right]^T $, so that
     \begin{equation*}
     df =  2\left(  \frac{\partial f}{\partial \mathbf{q}}\right)^T d\mathbf{q} + 2(d\mathbf{q}^*)^T\frac{\partial f}{\partial \mathbf{q}^*}  = 4\Re\left\{ \left( \frac{\partial f}{\partial \mathbf{q}}\right)^T d\mathbf{q} \right\}
     \end{equation*}
     Upon substituting  $\frac{\partial f}{\partial q^*} = \left(\frac{\partial f}{\partial q}\right)^*$ to give
     $     df = 4\Re\left\{ \left( \frac{\partial f}{\partial \mathbf{q}^*}\right)^H d\mathbf{q} \right\}     $
     and applying the Schwarz inequality\footnote{ The inner product of two vectors $\mathbf{a}$ and $\mathbf{b}$ satisfies
     $
     | \langle \mathbf{a},\mathbf{b}  \rangle| \leq \| \mathbf{a}\| \|\mathbf{ b}\|
     $,
      where the equality holds  when $\mathbf{a}$ is collinear with $\mathbf{b}$, that is, $\mathbf{a}= \lambda \mathbf{b}$.} we have
     \begin{eqnarray}
     df = 4 \Re | \langle  \frac{\partial f}{\partial \mathbf{q}^*},d\mathbf{q}    \rangle | \rightarrow
     df \leq 4  \|  \frac{\partial f}{\partial \mathbf{q}^*} \| \| d\mathbf{q} \|
     \label{eq:gradient:schwarz}
     \end{eqnarray}

\noindent \textbf{Remark\#8:}
     The maximum change of the differential $df$ occurs when $d\mathbf{q}$ is in the direction of $\frac{\partial f}{\partial \mathbf{q}^*}$, thus making the conjugate gradient $\nabla_{w^*}$ a natural choice of gradient in the optimization of real valued functions of quaternion variables.

\subsection{ Unifying the Quaternion Gradient: The I-Gradient}
As shown above, the gradient based on the $\mathbb{HR}$-derivatives is rigorous and has provided a quantum step forward in the derivation of stochastic gradient algorithms in the quaternion domain. However, despite its mathematical correctness it does not necessarily equip the quaternion domain with algorithms that are generic extensions of its real- and complex- valued counterparts, the LMS or CLMS (as seen from \cite{Clive_QLMS_2009} and (\ref{Eq:QLMS:w(k+1)Final})). It is therefore necessary to introduce a novel gradient definition that not only unifies the QLMS algorithms but also simplifies into CLMS and LMS for lower dimensional data. To that end, consider  the involution based representation in (\ref{eq:involution}), where \( q= \frac{1}{2}\left(q ^{{\imath}^*} + q ^{{\jmath}^*} + q ^{{\kappa}^*} -q^* \right)  \) and \( q^*= \frac{1}{2}\left(q^{\imath} + q^{\jmath} + q^{\kappa} -q \right)  \). In this way, we can  write the partial derivative \( \frac{\partial f}{\partial q}     \) and its conjugate $\left(\frac{\partial f}{\partial q}\right)^*$ as
\begin{equation}
 \frac{\partial f}{\partial q} =   \frac{1}{2} \left(  \left(\frac{\partial f}{\partial q}\right) ^{{\imath}^*} + \left(\frac{\partial f}{\partial q}\right) ^{{\jmath}^*} + \left(\frac{\partial f}{\partial q}\right) ^{{\kappa}^*} - \left(\frac{\partial f}{\partial q}\right)^*   \right)
 \label{Eq:SimpleQLMS:df/dq}
 \end{equation}
\begin{equation}
 \left(\frac{\partial f}{\partial q}\right)^* =   \frac{1}{2} \left(  \left(\frac{\partial f}{\partial q}\right)^{\imath} + \left(\frac{\partial f}{\partial q}\right)^{\jmath} + \left(\frac{\partial f}{\partial q}\right)^{\kappa} - \frac{\partial f}{\partial q}   \right)
 \label{Eq:SimpleQLMS:(df/dq)*}
 \end{equation}
 For a real valued function $f(q)$, using the identities in (\ref{Eq:QuaternionGradient:Identity})  and the fact that for a real function  of quaternion variables $  \frac{ \partial f}{\partial q}  = \left(\frac {\partial f(q)}{\partial q^*} \right)^* $, we arrive at
\begin{equation}
 \frac{\partial f}{\partial q} =   \frac{1}{2} \left(  \frac{\partial f}{\partial q ^{{\imath}^*}} + \frac{\partial f}{\partial q ^{{\jmath}^*}} + \frac{\partial f}{\partial q ^{{\kappa}^*}} - \frac{\partial f}{\partial q^*}   \right) \label{Eq:SimpleQLMS:(df/dq)2}
 \end{equation}
\begin{equation}
 \frac{\partial f}{\partial q^*} =   \frac{1}{2} \left(  \frac{\partial f}{\partial q^{\imath}} + \frac{\partial f}{\partial q^{\jmath}} + \frac{\partial f}{\partial q^{\kappa}} - \frac{\partial f}{\partial q}   \right)
 \label{Eq:SimpleQLMS:(df/dq)*2}
 \end{equation}
Observe from (\ref{eq:gradient:schwarz}), that it is the conjugate gradient
 \begin{equation*}
 \frac{\partial f}{\partial q^*} =   \frac{1}{2} \left(  \frac{\partial f}{\partial q^{\imath}} + \frac{\partial f}{\partial q^{\jmath}} + \frac{\partial f}{\partial q^{\kappa}} - \frac{1}{2} \left(  \frac{\partial f}{\partial q ^{{\imath}*}} + \frac{\partial f}{\partial q ^{{\jmath}*}} + \frac{\partial f}{\partial q ^{{\kappa}*}} - \frac{\partial f}{\partial q^*}   \right)   \right)
 \end{equation*}
that provides the maximum rate of change (steepest direction). The conjugate gradient can also be written as
  \begin{equation*}
 \frac{3}{4}\frac{\partial f}{\partial q*} =   \frac{1}{2} \left( \frac{1}{2}\left( \frac{\partial f}{\partial q^{\imath}} + \frac{\partial f}{\partial q^{\jmath}} + \frac{\partial f}{\partial q^{\kappa}}    \right) + \frac{1}{2}\left( \frac{\partial f}{\partial q^{\imath}}- \frac{\partial f}{\partial q ^{{\imath}*}} \right) +  \frac{1}{2}\left( \frac{\partial f}{\partial q^{\jmath}}- \frac{\partial f}{\partial q ^{{\jmath}*}} \right) + \frac{1}{2}\left( \frac{\partial f}{\partial q^{\kappa}}- \frac{\partial f}{\partial q ^{{\kappa}*}} \right)\right)
 \end{equation*}
and after substituting   $\frac{1}{2}(q - q^*)=\Im(q)$, where $\Im$ stands for the imaginary (vector) part of $q$, and using
    $\frac{\partial f(q)}{\partial q^*} = \left( \frac{\partial f(q)}{\partial q} \right)^*$, we have
   \begin{equation}
 \frac{\partial f}{\partial q^*} =   \frac{1}{3} \left( \left( \frac{\partial f}{\partial q^{\imath}} + \frac{\partial f}{\partial q^{\jmath}} + \frac{\partial f}{\partial q^{\kappa}}    \right) + 2 \left( \Im\left( \frac{\partial f}{\partial q^{\imath}}  \right) + \Im\left( \frac{\partial f}{\partial q^{\jmath}}  \right) + \Im\left( \frac{\partial f}{\partial q^{\kappa}}  \right) \right) \right)
 \label{Eq:SimpleQLMS:(df/dq)3}
 \end{equation}
 Notice that now both the real and imaginary parts are a function of only $ \frac{\partial f}{\partial q^{\imath}} $, $ \frac{\partial f}{\partial q^{\jmath}} $ and $ \frac{\partial f}{\partial q^{\kappa}} $, that is
 \begin{equation}
 \Re \left[\frac{\partial f}{\partial q^*} \right] + \Im \left[\frac{\partial f}{\partial q^*} \right]= \frac{1}{3}\Re \left[\frac{\partial f}{\partial q^{\imath}} + \frac{\partial f}{\partial q^{\jmath}} + \frac{\partial f}{\partial q^{\kappa}}\right] + \Im \left[\frac{\partial f}{\partial q^{\imath}} + \frac{\partial f}{\partial q^{\jmath}} + \frac{\partial f}{\partial q^{\kappa}}\right]
 \end{equation}
  and that the direction of the vector part of the gradient $\frac{\partial f}{\partial q^*}$ is equivalent to that within the $\mathbb{HR}$-calculus, that is, $\frac{\partial f}{\partial q^{\imath}} + \frac{\partial f}{\partial q^{\jmath}} + \frac{\partial f}{\partial q^{\kappa}}$. In this way, we introduce  a novel quaternion gradient, termed the involution- or I-gradient, given by
 \begin{equation}
 \nabla_q J = \frac{\partial J}{\partial q^{\imath}} + \frac{\partial J}{\partial q^{\jmath}} + \frac{\partial J}{\partial q^{\kappa}}
 \end{equation}
  Using (\ref{Eq:SimpleQLMS:(df/dq)*2}) and the  $\mathbb{HR}$- and $\mathbb{HR}^*$-derivatives in (\ref{Eq:HRcalculus:HRderivatives}) and (\ref{Eq:HRcalculus:HR*derivatives}), the I-gradient can be written in terms of the conjugate gradient as
  \begin{equation}
 \nabla_q J = \frac{\partial J}{\partial q^*} + \frac{1}{2}\frac{\partial J}{\partial q_r}
 \end{equation}
Notice that the I-gradient is equivalent to the $\mathbb{HR}$-gradient for 3D processes (pure quaternions), a typical  case in most real world applications.

\noindent \textbf{Remark\#9:} For 4D processes (full quaternion), the I-gradient includes an additional term $\frac{1}{2}\frac{\partial J}{\partial q_r}$ as compared to the $\mathbb{HR}$-gradient, thus exhibiting increased  steepness in the direction of the real component, and providing a faster convergence of the so derived learning algorithms. However, as the sign of the real gradient component remains unchanged, both gradients converge to the same optimum solution.

\section{A unifying treatment of quaternion valued adaptive algorithms \label{Sec:Unified_approach}}
As mentioned earlier, a major stumbling block in the development of quaternion-valued adaptive filtering algorithms has been the issue of quaternion gradient, which in its standard form introduces uncertainties in the position of the imaginary units and is thus not unique. We shall now illustrate that by using the quaternion conjugate gradient defined by $\nabla_{w^*} J$ (as given by the $\mathbb{HR}$-calculus),  a real function of quaternion variables can be differentiated directly, without resorting to tedious component-wise partial gradients or introducing any uncertainty in the placement of the unit vectors.

This also gives us the opportunity to unify a number of existing adaptive filtering algorithms in $\mathbb{H}$, and to provide a rigorous platform for future developments. In addition, we show how the I-gradient $\nabla_w J= \nabla_{w^{\imath}} J+ \nabla_{w^{\jmath}} J + \nabla_{w^{\kappa}} J$, derived in Section \ref{Sec:QuaternionGradient}  can be used to introduce generic quaternion valued extensions of the existing real and complex adaptive filters, both strictly linear and widely linear ones. To be consistent with the original derivation of the QLMS in \cite{Clive_QLMS_2009}, the filter output assumes the form $y(k)=\mathbf{w}^T(k) \mathbf{x}(k)$, however, any other inner product of $\mathbf{w}(k)$ and $\mathbf{x}(k)$, e.g. $\mathbf{x}^T(k)\mathbf{w}(k)$, can be equally used to derive the algorithms.

\subsection{The Wiener Solution}
The Wiener filter aims to find the optimal coefficient (weight) vector that minimizes the mean square error, or equivalently to reach the minimum of the cost function, given by
\vspace{0mm}
\begin{equation}
J(k) = E[e(k)e^*(k)] = E[|e(k)|^2]
\label{Eq:Wiener:cost function}
\vspace{0mm}
\end{equation}
where $e(k)$ is the output error of an adaptive filter
\vspace{0mm}
\begin{equation*}
e(k)= d(k) - \mathbf{w}^T \mathbf{x}(k)
\vspace{0mm}
\end{equation*}
and $\mathbf{w} \in \mathbb{H}^{N \times 1}$ and $d(k)  \in \mathbb{H}$ are respectively  the filter weight vector and teaching signal.
To obtain the optimal weight vector $\mathbf{w}_o$, we need to calculate the gradient (using the conjugate gradient in (\ref{Eq:HRcalculus:HR*derivatives})) in the form\footnote{To simplify the calculation of the quaternion gradient, the product rule is used. See Appendix \ref{app:product_rule} for a justification for using the product rule.}
\begin{equation}
 \nabla_{\mathbf{w}^*}J= \frac{\partial J(k)}{\partial \mathbf{w}^*}= E\left[e(k)\frac{\partial e^*(k)}{\partial \mathbf{w}^*} + \frac{\partial e(k)}{\partial \mathbf{w}^*}e^*(k)\right]
 \label{eq:wiener:1}
\end{equation}
Normally, this is not straightforward to achieve using the pseudogradient, however, by virtue of the $\mathbb{HR}$-calculus in  (\ref{eq:wiener:1}), we directly obtain
\begin{equation*}
\frac{\partial J(k)}{\partial \mathbf{w}^*}=    E\left[\frac{1}{2}\mathbf{x}(k)e^*(k) -  e(k)\mathbf{x}^*(k)\right]
\end{equation*}
Substitute for $e(k)$ and set to zero to give
\begin{equation*}
E\left[\frac{1}{2}  \mathbf{x}(k) \big(   d^*(k) -  \mathbf{x}^H(k)\mathbf{w}^*   \big) -  \big(   d(k) - \mathbf{w}^T\mathbf{x}(k)   \big)\mathbf{x}^*(k) \right]= \mathbf{0}
\end{equation*}
Upon rearranging the terms above we have
\begin{eqnarray*}
E\big[  \mathbf{w}^T \mathbf{x}(k)\mathbf{x}^*(k)   -   \frac{1}{2} \mathbf{x}(k) \mathbf{x}^H(k)\mathbf{w}^*  \big] =
E\big[d(k)\mathbf{x}^*(k) - \frac{1}{2}  \mathbf{x}(k)d^*(k) \big]
\end{eqnarray*}
or in a more compact form
\begin{equation}
\left(\mathbf{w}^T\mathbf{R}_{x}\right)^T -  \frac{1}{2}\mathbf{R}_{x}\mathbf{w}^* = \mathbf{r}_{dx^*} - \frac{1}{2}\mathbf{r}_{xd^*}
\label{Eq:RLS:wRxx}
\end{equation}
Noting that $\left(\mathbf{w}^T\mathbf{R}_{x}\right)^H = \mathbf{R}_{x}\mathbf{w}^*$ and $(\mathbf{r}_{dx^*})^* = \mathbf{r}_{xd^*}$, we arrive at the strictly linear \emph{quaternion Wiener filter} in the form
\vspace{0mm}
\begin{equation}
\mathbf{w}^T_o = \mathbf{r}^T_{dx^*}\mathbf{R}^{-1}_{x}
\label{Eq:RLS:wiener}
\end{equation}
\vspace{0mm}
Notice the subtle difference from the real- and complex-valued Wiener solution - the quaternion solution is given in terms of the transpose of $\mathbf{w}$ as the rigid noncommutativity of $\mathbf{r}^T_{dx^*}\mathbf{R}^{-1}_{x}$ means that $\mathbf{R}^{-1}_{x} \mathbf{r}_{dx^*}$ would represent a different solution.

\subsection{The Widely Linear Wiener Solution}
When the system model is widely linear as in (\ref{eq:WL-model}), the augmented filter weight takes the form $\mathbf{w}^a=[\mathbf{u}^T,\mathbf{v}^T,\mathbf{g}^T,\mathbf{h}^T]^T$ and the augmented input vector $\mathbf{x}^a=[\mathbf{x}^{T},\mathbf{x} ^{{\imath}T},\mathbf{x} ^{{\jmath}T},\mathbf{x} ^{{\kappa}T}]$. Following on (\ref{Eq:RLS:wiener}), the \emph{widely linear Wiener} solution then becomes
\begin{equation}
\mathbf{w}^{aT}_o= \mathbf{r}_{dx^*}^T\left(\mathbf{R}^{a}_{x}\right)^{-1}
\end{equation}
 This model represents an optimal second order solution for both proper and improper signals, whereas the standard, strictly linear, Wiener filter in (\ref{Eq:RLS:wiener}) is optimal only for proper signals, for which the pseudocovariances matrices $\mathbf{P}_x$, $\mathbf{S}_x$ and $\mathbf{T}_x$ in (\ref{eq:widely_linear_Rqq}) vanish.

\subsection{The Quaternion Least Mean Square Algorithm}
We now illustrate the advantages  of the quaternion gradients introduced in this work, by providing an elegant derivation of the strictly linear QLMS, for which the original derivation  \cite{Clive_QLMS_2009} used cumbersome component-wise pseudogradients. Starting from the cost function in (\ref{Eq:Wiener:cost function}), the stochastic gradient weight update is given by ($\mu$ is the real-valued step size)
\vspace{-2mm}
\begin{equation}
\textbf{w}(k+1) = \textbf{w}(k) - \mu\nabla_\textbf{w} J(k)
\label{Eq:QLMS:w(k+1)}
\end{equation}
\vspace{-2mm}and noting that
\vspace{-2mm}
\begin{eqnarray*}
e(k)=d(k) - \textbf{w}^T(k)\textbf{x}(k)  \hspace{20mm}
e^*(k)=d^*(k) - \textbf{x}^H(k)\textbf{w}^*(k)
\end{eqnarray*}
\vspace{-2mm}the gradient of the cost function becomes
\begin{equation}
\nabla_{\textbf{w}^*} J(k)= \frac{1}{2} \left( e(k) \frac{\partial e^*(k)}{\partial \textbf{w}^*(k)} + \frac{\partial e(k)}{\partial \textbf{w}^*(k)}e^*(k) \right)
\label{Eq:QLMS:J}
 \end{equation}
Using the $\mathbb{HR}$-gradient $\nabla_{\textbf{w}^*}$, we can differentiate \(e(k)\) and \(e^*(k)\) directly, without resorting to the partial derivatives with respect to $\mathbf{w}_r$, $\mathbf{w}_{\imath}$, $\mathbf{w}_{\jmath}$ and $\mathbf{w}_{\kappa}$, thus giving
 \begin{equation}
 \frac{\partial e^*(k)}{\partial \textbf{w}^*} = -\textbf{x}^*(k)
 \label{Eq:QLMS:de*/dw*}
 \end{equation}
On the other hand, to calculate \(\frac{\partial e(k)}{\partial \textbf{w}^*(k)}  \), we must first write \(\textbf{w}\) in terms of the conjugate involutions \(\textbf{w}^*,\textbf{w} ^{{\imath}*},\textbf{w} ^{{\jmath}*},\textbf{w} ^{{\kappa}*} \), so that the argument takes the form required by (\ref{Eq:HRcalculus:HR*derivatives}), given by
\begin{equation}
e(k) = d(k) - \frac{1}{2}\big({\textbf{w} ^{{\imath}*}}(k) + {\textbf{w} ^{{\jmath}*}}(k) + {\textbf{w} ^{{\kappa}*}}(k) - \textbf{w}^*(k)\big)^T\textbf{x}(k)
\end{equation}
Using the $\mathbb{HR}$-calculus in (\ref{Eq:HRcalculus:HR*derivatives}) to perform the direct differentiation of $e(k)$ gives
\begin{equation}
\frac{\partial e(k)}{\partial \textbf{w}^*} = \frac{1}{2}\textbf{x}(k)
\label{Eq:QLMS:de/dw*}
\end{equation}
and combining the terms above yields the update, referred to as the HR-QLMS, in the form
\begin{equation}
\textbf{w}(k+1) = \textbf{w}(k) + \mu\left( \frac{1}{2}e(k)\textbf{x}^*(k) - \frac{1}{4}\textbf{x}(k)e^*(k) \right)
\label{Eq:QLMS:w(k+1)Final}
\end{equation}
A comparison with the original QLMS algorithm in \cite{Clive_QLMS_2009} whose form is
\vspace{0mm}
\begin{equation}
\textbf{w}(k+1) = \textbf{w}(k) + \mu\big( 2e(k)\textbf{x}^*(k) - \textbf{x}^*(k)e^*(k) \big)
\label{Eq:QLMS:w(k+1)standard}
\vspace{0mm}
\end{equation}
leads to the following observations.

\noindent \textbf{Remark\#10:}
 The standard QLMS weight update term is four times larger than that within the HR-QLMS, the difference arising due to the component-wise gradient calculation in the standard QLMS. However, this difference can be absorbed into the learning rate $\mu$.

\noindent \textbf{Remark\#11:}
 The second term in the original QLMS weight update uses  \(\textbf{x}^*(k)\) instead of \(\textbf{x}(k)\); this difference arises due to the interpretation of the gradient \(\frac{\partial f(\textbf{w})}{\partial \textbf{w}^*}\) in (\ref{Eq:HRcalculus:HR*derivatives}). In the standard QLMS in \cite{Clive_QLMS_2009} it is assumed that the unit vectors \(\imath\), \(\jmath\) and \(\kappa\) are always on the right hand side of the real component-wise derivatives $\frac{\partial f}{\mathbf{q}_r}$, $\frac{\partial f}{\mathbf{q}_{\imath}}$, $\frac{\partial f}{\mathbf{q}_{\jmath}}$ and $\frac{\partial f}{\mathbf{q}_{\kappa}}$. In the QLMS derivation based on the $\mathbb{HR}$-gradient,  having the imaginary units before the derivatives $\frac{\partial f}{\mathbf{q}_r}$, $\frac{\partial f}{\mathbf{q}_{\imath}}$, $\frac{\partial f}{\mathbf{q}_{\jmath}}$ and $\frac{\partial f}{\mathbf{q}_{\kappa}}$ for the gradient $\frac{\partial e}{\partial \mathbf{w}^*}$ and after  the partial derivatives for the gradient $\frac{\partial e^*}{\partial \mathbf{w}^*}$ results in the HR-QLMS weight update in (\ref{Eq:QLMS:w(k+1)Final}).

 Section \ref{Section:equivalence} provides an in-depth comparative analysis, and shows that these two expressions effectively provide the same solution.

\subsection{The Widely Linear Least Mean Square Algorithm (WL-QLMS)}
Using the widely linear model in Section \ref{Sec:WidelyLinearModel} and following on  the analysis of the QLMS class of algorithms, we can now rederive the WL-QLMS \cite{Clive_WL-QLMS_2009} in an efficient way, and establish an elegant and rigorous framework for the derivation of widely linear quaternion valued algorithms that account for the improperness of quaternion data. We assume that the output of the filter takes the form
$
y(k)= \textbf{w}^{aT}(k)\textbf{x}^a(k)
$
where the augmented weight vector $\textbf{w}^a(k)=[\textbf{u}^T(k) ,  \textbf{v}^T(k)  ,  \textbf{g}^T(k)  ,  \textbf{h}^T(k) ]^T$  and the augmented input vector $\mathbf{x}^a(k)= [\mathbf{x}^T(k)  ,  \mathbf{x} ^{{\imath}T}(k)  ,  \mathbf{x} ^{{\jmath}T}(k)  ,  \mathbf{x} ^{{\kappa}T}(k)]^T$. The weight updates are  calculated from (\ref{Eq:QLMS:w(k+1)}), based on the gradient
 \begin{eqnarray}
\nabla_\textbf{w} J(k) &=& \frac{1}{2} \left( e(k)  \frac{\partial e^*(k)}{\partial \textbf{w}^{a*}(k)} + \frac{\partial e(k)}{\partial \textbf{w}^{a*}(k)}e^*(k) \right)
\label{Eq:WL-QLMS:dJ/duvgh}
\end{eqnarray}
By virtue of the $\mathbb{HR}$-calculus, the partial derivatives $\left[ \frac{\partial e^*(k)}{\partial \textbf{u}^*(k)}, \frac{\partial e^*(k)}{\partial \textbf{v}^*(k)}, \frac{\partial e^*(k)}{\partial \textbf{g}^*(k)},  \frac{\partial e^*(k)}{\partial \textbf{h}^*(k)} \right]^T $  within the gradient $\frac{\partial e^*(k)}{\partial \textbf{w}^{a*}(k)}$ can be calculated by direct differentiation, giving
 \begin{eqnarray}
 \frac{\partial e^*(k)}{\partial \textbf{u}^*(k)} = -\textbf{x}^*(k) \hspace{10mm} \frac{\partial e^*(k)}{\partial \textbf{v}^*(k)} = -\textbf{x} ^{{\imath}^*}(k) \hspace{10mm} 
 \frac{\partial e^*(k)}{\partial \textbf{g}^*(k)} = -\textbf{x} ^{{\jmath}^*}(k) \hspace{10mm} \frac{\partial e^*(k)}{\partial \textbf{h}^*(k)} = -\textbf{x} ^{{\kappa}^*}(k)
 \label{Eq:WL-QLMS:de*/duvgh}
 \end{eqnarray}
To obtain the partial derivatives $ \frac{\partial e(k)}{\partial \textbf{w}^{a*}(k)} = \left[ \frac{\partial e(k)}{\partial \textbf{u}^*(k)}, \frac{\partial e(k)}{\partial \textbf{v}^*(k)}, \frac{\partial e(k)}{\partial \textbf{g}^*(k)} \frac{\partial e(k)}{\partial \textbf{h}^*(k)} \right]$, substitute \\ $\mathbf{w}^a=\frac{1}{2}\left((\mathbf{w}^a) ^{{\imath}*} + (\mathbf{w}^a) ^{{\jmath}*} + (\mathbf{w}^a) ^{{\kappa}*} - \mathbf{w}^{a*} \right)$, to yield
\begin{eqnarray}
 \frac{\partial e(k)}{\partial \textbf{u}^*(k)} = \frac{1}{2}\textbf{x}(k) \hspace{10mm}
 \frac{\partial e(k)}{\partial \textbf{v}^*(k)}= \frac{1}{2}\textbf{x}^{\imath}(k) \hspace{10mm} 
 \frac{\partial e(k)}{\partial \textbf{g}^*(k)} = \frac{1}{2}\textbf{x}^{\jmath}(k) \hspace{10mm}
 \frac{\partial e(k)}{\partial \textbf{h}^*(k)}= \frac{1}{2}\textbf{x}^{\kappa}(k)
 \label{Eq:WL-QLMS:de/duvgh}
 \end{eqnarray}
 finally arriving at the weight update for the WL-QLMS algorithm in the form
\begin{eqnarray}
\textbf{w}^a(k+1) &=& \textbf{w}^a(k) + \mu\left(\frac{1}{2}e(k)\textbf{x}^{a*}(k) -  \frac{1}{4}\textbf{x}^a(k)e^*(k) \right)
\label{eq:WLQLMS:weightupdate}
\end{eqnarray}
This expression for the WL-QLMS differs from the existing WL-QLMS algorithm introduced in \cite{Clive_WL-QLMS_2009}, the differences arising for the same reason as those for the QLMS.

\subsection{ A Class of Generic QLMS Algorithms \label{SEc:IQLMS}}

Upon applying the I-gradient derived in Section \ref{Sec:QuaternionGradient} to  the cost function in (\ref{Eq:Wiener:cost function}), we obtain
\begin{eqnarray}
\nabla_\textbf{w} J(k) &=& \frac{\partial J}{\partial w^{\imath}} + \frac{\partial J}{\partial w^{\jmath}} + \frac{\partial J}{\partial w^{\kappa}} \nonumber \\
&=& \frac{1}{2} \left( e(k) \frac{\partial e^*(k)}{\partial \textbf{w}^{\imath}(k)} + \frac{\partial e(k)}{\partial \textbf{w}^{\imath}(k)}e^*(k) +  e(k) \frac{\partial e^*(k)}{\partial \textbf{w}^{\jmath}(k)} + \frac{\partial e(k)}{\partial \textbf{w}^{\jmath}(k)}e^*(k)  \right. \nonumber \\
&+& \left. e(k) \frac{\partial e^*(k)}{\partial \textbf{w}^{\kappa}(k)} + \frac{\partial e(k)}{\partial \textbf{w}^{\kappa}(k)}e^*(k)\right)
\label{Eq:SimpleQLMS:dJ/dwi}
 \end{eqnarray}
Using the $ \mathbb{HR}$-calculus,
$
\frac{\partial e(k)}{\partial \mathbf{w}^{\imath }} = \frac{\partial e(k)}{\partial \mathbf{w}^{\jmath }} = \frac{\partial e(k)}{\partial \mathbf{w}^{\kappa}} = \mathbf{0}
\label{Eq:SimpleQLMS:de/dwi}
$,
and similarly to the derivation of QLMS, we can write
\begin{equation*}
e^*(k) = d^*(k) - \frac{1}{2}\big(\textbf{w}^{\imath }(k) + \textbf{w}^{\jmath }(k) + \textbf{w}^{\kappa }(k) - \textbf{w}(k) \big)^T\textbf{x}(k)
\end{equation*}
By direct differentiation,
$\frac{\partial e^*(k)}{\partial \mathbf{w}^{\imath}}= \frac{\partial e^*(k)}{\partial \mathbf{w}^{\jmath}} = \frac{\partial e^*(k)}{\partial \mathbf{w}^{\kappa}} = -\frac{1}{2}\mathbf{x}(k)
\label{Eq:SimpleQLMS:de*/dwi}
$
yielding the weight update
\begin{equation}
\mathbf{w}(k+1)= \mathbf{w}(k) + \frac{3}{4}\mu e(k)\mathbf{x}^*(k)
\label{Eq:SimpleQLMS:CLMS_type}
\end{equation}
This algorithm  is termed the IQLMS (the I-gradient based QLMS).

\noindent\textbf{Remark\#12}: The IQLMS  has the same generic form as the LMS and CLMS \cite{Widrow_CLMS_1975}, as the term $\frac{3}{4}$ can be absorbed into the learning rate.

\noindent\textbf{Remark\#13}: The IQLMS requires half the computations of the QLMS.

In the same spirit, from (\ref{eq:WLQLMS:weightupdate}) in order to derive the WL-IQLMS, a generic extension of the complex widely linear ACLMS \cite{Danilo_ACLMS} algorithm, we employ the I-gradient to obtain
\begin{equation}
\mathbf{w}^a(k+1)= \mathbf{w}^a(k) + \frac{3}{4}\mu e(k)\mathbf{x}^{a*}(k)
\label{Eq:SimpleQLMS:CLMS_type}
\end{equation}
where $ \mathbf{x}^a$ is the augmented input vector. It is easy to show, following the approach in \cite{Benesty_Duality} (see Appendix \ref{app:duality}) that the WL-IQLMS is structurally identical to the four-channel real LMS, whereas the WL-QLMS is not.

\noindent\textbf{Remark\#14}: Note  that the Wiener filter is already in its generic form and hence the corresponding algorithm derived using the I-gradient would be the same.

\section{Equivalence of the various forms of QLMS  \label{Section:equivalence}}
 Starting from the expression for the original QLMS  given by \cite{Clive_QLMS_2009}
\begin{equation}
\mathbf{w}(k+1)=\mathbf{w}(k) + \mu\left( \frac{1}{2}e(k)\mathbf{x}^*(k) - \frac{1}{4}\mathbf{x}^*(k)e^*(k)  \right)
\label{eq:equivalence1}
\end{equation}
together with the expression for the proposed HR-QLMS
\begin{equation}
\mathbf{w}(k+1)=\mathbf{w}(k) + \mu\left( \frac{1}{2}e(k)\mathbf{x}^*(k) - \frac{1}{4}\mathbf{x}(k)e^*(k)  \right)
\label{eq:equivalence:HR-QLMS}
\end{equation}
and the IQLMS algorithm
\begin{equation}
\mathbf{w}(k+1)= \mathbf{w}(k) + \frac{3}{4}\mu e(k) \mathbf{x}^*(k)
\label{eq:equivalence2}
\end{equation}
we can now rewrite the QLMS algorithms (\ref{eq:equivalence1})-(\ref{eq:equivalence2}) in a form that allows for a convenient analytical comparison.

Based on (\ref{eq:equivalence1}), from $\mathbf{x}^*(k) = \mathbf{x}(k) - 2\Im[\mathbf{x}(k)]$ and by splitting the term $\frac{1}{2}e(k)\mathbf{x}^*(k)$ we have
\begin{eqnarray*}
\mathbf{w}(k+1)&=&\mathbf{w}(k) + \frac{1}{4}\mu e(k)\mathbf{x}^*(k) + \frac{1}{4}\mu e(k)\mathbf{x}^*(k) - \frac{1}{4}\mu\mathbf{x}(k)e^*(k)  + \frac{1}{2}\mu \Im[\mathbf{x}(k)]e^*(k) \\
&=& \mathbf{w}(k) + \mu\frac{1}{2}\Im[ e(k)\mathbf{x}^*(k)] + \frac{1}{4}\mu e(k)\mathbf{x}^*(k) + \frac{1}{2}\mu \Im[\mathbf{x}(k)]e^*(k) \\
&=& \mathbf{w}(k) + \mu\frac{1}{2}\Im[ e(k)\mathbf{x}^*(k)] + \frac{1}{4}\mu e(k)\mathbf{x}^*(k) - \frac{1}{2}\mu\big[ e(k)\Im[\mathbf{x}(k)]\big]^*
\end{eqnarray*}
The real part of the weight update can therefore be  evaluated as
\begin{equation*}
 \Re[\mathbf{w}(k+1)] = \Re[\mathbf{w}(k)] + \frac{1}{4}\mu \Re[e(k)\mathbf{x}^*(k)] - \frac{1}{2}\mu \Re\big[ e(k)\Im[\mathbf{x}(k)]\big]
\end{equation*}
which, upon using $\mathbf{x}(k) = \Re[\mathbf{x}(k)] + \Im[\mathbf{x}(k)]$, yields
\begin{equation}
 \Re[\mathbf{w}(k+1)] = \Re[\mathbf{w}(k)] + \frac{1}{4}\mu \Re \big[e(k)\Re[\mathbf{x}(k)]\big] - \frac{3}{4}\mu \Re\big[ e(k)\Im[\mathbf{x}(k)]\big]
\end{equation}
Similarly, for the imaginary part of the QLMS weight update in (\ref{eq:equivalence1}), we have
\begin{equation*}
 \Im[\mathbf{w}(k+1)] = \Im[\mathbf{w}(k)] + \frac{3}{4}\mu \Im[e(k)\mathbf{x}^*(k)] + \frac{1}{2}\mu \Im\big[ e(k) \Im[\mathbf{x}(k)]\big]
\end{equation*}
\begin{equation}
 \Im[\mathbf{w}(k+1)] = \Im[\mathbf{w}(k)] + \frac{3}{4}\mu \Im\big[e(k)\Re[\mathbf{x}(k)]\big] - \frac{1}{4}\mu\big[ e(k)\Im[\mathbf{x}(k)]\big]
\end{equation}

Repeating the above process  for the weight updates in (\ref{eq:equivalence:HR-QLMS}) and (\ref{eq:equivalence2}), it can be shown that all the various forms of the  updates are based on a combination of the terms, $ + \left[e(k)\Re[\mathbf{x}(k)] \right]$ and $-\left[e(k)\Im[\mathbf{x}(k)] \right]$, or more explicitly
\begin{eqnarray}
\Re[\mathbf{w}(k+1)] = \Re[\mathbf{w}(k)] + a \text{ }\mu\Re\left[e(k)\Re[\mathbf{x}(k)] \right] - b\text{ } \mu \Re\left[e(k)\Im[\mathbf{x}(k)] \right]
\label{eq:equivalence3}
\end{eqnarray}
\begin{eqnarray}
\Im[\mathbf{w}(k+1)] =  \Im[\mathbf{w}(k)] + c\text{ } \mu\Im\left[e(k)\Re[\mathbf{x}(k)] \right] - d\text{ } \mu \Im\left[e(k)\Im[\mathbf{x}(k)] \right]
\label{eq:equivalence4}
\end{eqnarray}
Therefore, the weight updates for QLMS, HR-QLMS and IQLMS  differ only in the weighting factors $a$, $b$, $c$ and $d$, that is
\begin{itemize}
  \item For the QLMS derived in \cite{Clive_QLMS_2009}, the weighting coefficients are  $a=d=\frac{1}{4}$ and $b=c=\frac{3}{4}$.
  \item  For the QLMS in (\ref{eq:equivalence:HR-QLMS}), the weighting coefficients are $a=b=\frac{1}{4}$ and $c=d=\frac{3}{4}$.
  \item  For the IQLMS in (\ref{eq:equivalence2}), the weighting coefficients are $a=b=c=d=\frac{3}{4}$.
\end{itemize}
\noindent\textbf{Remark\#15}: The three QLMS algorithms in (\ref{eq:equivalence1})-(\ref{eq:equivalence2}) represent the same solution,  as the weight updates differ only in  the different weightings applied to the four components of the weight vector. This also suggests that the higher complexity of QLMS as opposed to IQLMS does not automatically result in a better performance.

\noindent\textbf{Remark\#16}: The different weightings in the weight update terms  mean that although all the algorithms will arrive at the same Wiener solution, the three versions of the QLMS will  exhibit differences in their convergence and steady state performances. For instance, we would expect the IQLMS to have a faster initial convergence than the  QLMS, due to the larger weighting of the weight vector components.

 Indeed, the analysis in Section \ref{Sec:Convergence} and simulations in Section \ref{sec:simulations} show that all the three algorithms have similar performances, yet the IQLMS requires only half the operations.
\section{Convergence analysis  \label{Sec:Convergence}}

We now provide a unified platform for the  convergence analysis of the class of QLMS algorithms. For rigour, this is achieved starting from  the stability conditions through to the mean square error performance for noncircular data of both strictly linear and widely linear filters. The mean square error analysis is conducted for the IQLMS and also provides a bound on the performance of the QLMS and HR-QLMS.

\subsection{Convergence in the Mean of the Strictly Linear IQLMS Algorithm}
To simplify the derivation, the output is re-written in  the form $\mathbf{w}^H\mathbf{x}$; this does not change any of its properties, but makes the analysis more mathematically tractable.
Without loss of generality, assume first that the teaching signal for the IQLMS is given by (output of the strictly linear model)
\vspace{-2mm}
\begin{equation*}
d(k)=\mathbf{w}^H_o\mathbf{x}(k) + n(k)
\end{equation*}
\vspace{-2mm}where $\mathbf{w}_o$ is the optimal filter weight vector and $n(k)$ quaternion quadruply white Gaussian noise\footnote{Each component of  quaternion noise is drawn from an independent Gaussian distribution with equal power.}.
The output error can therefore be written as
\vspace{-2mm}
\begin{equation}
e(k)=\mathbf{w}^H_o\mathbf{x}(k) + n(k) - \mathbf{w}^H(k)\mathbf{x}(k)
\label{eq:converg.analysis:e}
\vspace{-2mm}
\end{equation}
Substitute into the IQLMS weight update  in (\ref{eq:equivalence2}) to obtain
\begin{equation*}
\mathbf{w}(k+1)=\mathbf{w}(k) + \frac{3}{4}\mu \mathbf{x}(k)\mathbf{x}^H(k)\mathbf{w}_o(k) + \frac{3}{4}\mu \mathbf{x}(k)n^*(k) - \frac{3}{4}\mu \mathbf{x}(k)\mathbf{x}^H(k)\mathbf{w}(k)
\end{equation*}
and subtract $\mathbf{w}_o$ from both sides to give the weight error vector   $\mathbf{r}(k)=\mathbf{w}(k)-\mathbf{w}_o$ in the form
\begin{equation*}
\mathbf{r}(k+1)=\mathbf{r}(k) - \frac{3}{4}\mu \mathbf{x}(k)\mathbf{x}^H(k)\left(\mathbf{w}(k) - \mathbf{w}_o(k)\right) + \frac{3}{4}\mu \mathbf{x}(k)n^*(k)
\end{equation*}
For convergence in the mean, take the statistical expectation of both sides to yield
\begin{eqnarray}
E[\mathbf{r}(k+1)] &=& \left(\mathbf{I} - \frac{3}{4}\mu E[\mathbf{x}(k)\mathbf{x}^H(k)]\right)E[\mathbf{r}(k)] + \frac{3}{4}E[\mathbf{x}(k)n^*(k)]
= (\mathbf{I} - \frac{3}{4}\mu \mathbf{R}_{x})E[\mathbf{r}(k)]
\label{eq:IQLMS_error_update}
\end{eqnarray}
The recursion for the weight error vector $\mathbf{r}(k)$ converges  for $||\mathbf{I} - \frac{3}{4}\mu \mathbf{R}_{x}||<1$, where the norm is defined as the spectral norm (corresponding to the largest absolute eigenvalue).
Making use of the unitary transform $\mathbf{R}_x=\mathbf{Q}\boldsymbol{\Lambda} \mathbf{Q}^H$, where $\mathbf{Q}$ is the matrix of eigenvectors and $\boldsymbol{\Lambda}$ is the eigenvalue matrix, and  rotating $\mathbf{r}(k)$ to obtain  $\mathbf{r}'(k)=\mathbf{Q}^HE[\mathbf{v}(k)]$ yields
\vspace{-2mm}
\begin{equation*}
\mathbf{r'}(k+1) = (\mathbf{I} - \frac{3}{4}\mu\boldsymbol{\Lambda})\mathbf{r'}(k)
\vspace{-2mm}
\end{equation*}
Since $\mathbf{I} - \frac{3}{4}\mu \mathbf{\Lambda}  $ is diagonal, the  `modes of convergence' can be expressed as
\begin{equation}
r'_n(k+1) =   \left(1 - \frac{3}{4}\mu \lambda_n(\mathbf{R}_{x})  \right) r'_n(k) \hspace{5mm} r'_n \in \mathbf{r}, \text{ }n=1,...,N
\end{equation}
Therefore, in order for every $r_n(k) \in \mathbf{r}(k)$ to converge
\begin{equation*}
|1 - \frac{3}{4}\mu \lambda_{max}(\mathbf{R}_{x})|  <1
\end{equation*}
\begin{equation}
 0 < \mu <\frac{8}{3\lambda_{max}(\mathbf{R}_{x})}
 \label{Eq:IQLMS_Convergence:range}
\end{equation}
where the symbol $\lambda_{max}(\cdot)$ denotes the maximum eigenvalue.
Using the relationship between the trace of a matrix and its eigenvalues, that is, $Tr(\mathbf{R}) = \sum_{i=1}^N \lambda_i$, the upper bound on the learning rate can be conveniently approximated by
 \begin{equation*}
 0 < \mu <  \frac{8}{3Tr(\mathbf{R}_{x})}
\end{equation*}
where  $Tr(\cdot)$ denotes matrix trace operation.

\subsection{Convergence in the Mean of HR-QLMS}

We shall first rewrite the weight update of the QLMS  as
\begin{equation}
\mathbf{w}(k+1)= \mathbf{w}(k) + \mu\left(\frac{1}{2}e(k)\mathbf{x}^*(k) - \frac{1}{4}\mathbf{x}(k)e^*(k) \right)
\label{Eq:HRQLMS_Convergence:w(k+1)}
\end{equation}
By substituting for $e(k)$ into (\ref{Eq:HRQLMS_Convergence:w(k+1)}) and for the weight error vector  $\mathbf{r}(k)=\mathbf{w}(k)-\mathbf{w}_o(k)$, it can be shown (see Appendix \ref{app:proof_convergence_HR_v}) that the weight error vector takes the form
\begin{equation}
\mathbf{r}^T(k+1) = \mathbf{r}^T(k) -  \frac{1}{2} \mu\left(\mathbf{r}^T(k)\mathbf{x}(k)\mathbf{x}^H(k)\right) +  \frac{1}{4} \mu\left(\mathbf{r}^T(k)\mathbf{x}(k)\mathbf{x}^H(k)\right)^* +\frac{1}{2} \mu n(k) \mathbf{x}^H(k) - \frac{1}{4} \mu \mathbf{x}^T(k)n^*(k)
\label{eq:convergence:HRQLMS:v}
\end{equation}
It is not possible to obtain a recursive expression for $\mathbf{r}(k)$, as for the IQLMS,  due to the presence of the conjugate terms $\left(\mathbf{r}^T(k)\mathbf{x}(k)\mathbf{x}^H(k)\right)^*$. Instead, by splitting $\mathbf{r}(k)$ into its real and imaginary components (see Appendix \ref{App:convergence}) and taking the statistical expectation, the following recursive expression can be obtained
\begin{equation}
\boldsymbol{\omega}^T(k+1)= \boldsymbol{\omega}^T(k)[\mathbf{I} - \frac{3}{4}\mu\mathbf{R}_b]
\end{equation}
where $\boldsymbol{\omega}(k) = [E[\mathbf{r}^T_r(k)] \text{,} \imath E[\mathbf{r}^T_\imath(k)] \text{,} \jmath E[\mathbf{r}^T_\jmath(k)] \text{,} \kappa E[\mathbf{r}^T_\kappa(k)]]^T $
and
\[
\mathbf{R}_{b}=
\left( \begin{array}{cccc}
\frac{1}{3}[\mathbf{R}_{x}]_r & {\imath}[\mathbf{R}_{x}]_{\imath} & {\jmath}[\mathbf{R}_{x}]_{\jmath} & {\kappa}[\mathbf{R}_{x}]_{\kappa} \\
{\imath}\frac{1}{3}[\mathbf{R}_{x}]_{\imath} & [\mathbf{R}_{x}]_r  & {\kappa}[\mathbf{R}_{x}]_{\kappa} & {\jmath}[\mathbf{R}_{x}]_{\jmath} \\
{\jmath}\frac{1}{3}[\mathbf{R}_{x}]_{\jmath} & {\kappa}[\mathbf{R}_{x}]_{\kappa} & [\mathbf{R}_{x}]_r  & {\imath}[\mathbf{R}_{x}]_{\imath} \\
{\kappa}\frac{1}{3}[\mathbf{R}_{x}]_{\kappa} & {\jmath}[\mathbf{R}_{x}]_{\jmath} & {\imath}[\mathbf{R}_{x}]_{\imath} & [\mathbf{R}_{x}]_r \\
\end{array} \right)
\]
and the notation $[\mathbf{R}_{x}]_{\eta}$ indicates the $\eta$ part of $\mathbf{R}_{x}$ for $\eta=\{r,\imath,\jmath,\kappa\}$.
Similarly to the IQLMS, the condition for convergence in the mean becomes
\begin{equation}
|\mathbf{I} - \frac{3}{4}\mu\mathbf{R}_b|<1
\end{equation}
while the bound on the  step size in terms of the eigenvalues of $\mathbf{R}_b$, becomes
      \begin{equation}
 0 < \mu <\frac{8}{3\lambda_{max}(\mathbf{R}_b)}
\end{equation}
To compare the performance of the IQLMS and QLMS, the eigenvalues of the matrix $\mathbf{R}_b$ must be related to the eigenvalues of the matrix $\mathbf{R}_{x}$. As shown  in Appendix \ref{App:convergence_Proof}, this gives
\begin{eqnarray}
\lambda_{max}(\mathbf{R}_b) \approx \lambda_{max}(\mathbf{R}_{x}) \label{convergence:HRQLMS:max}\\
\lambda_{min}(\mathbf{R}_b) \approx \frac{1}{3}\lambda_{min}(\mathbf{R}_{x}) \label{convergence:HRQLMS:min}
\end{eqnarray}

\subsection{Convergence in the Mean of QLMS}

Starting from the weight update of the QLMS in (\ref{eq:equivalence1}) it can be shown (see Appendix \ref{app:proof_convergence:Clive_v}) that the weight error vector takes the form
\begin{equation}
\mathbf{r}^T(k+1) = \mathbf{r}^T(k) -  \frac{1}{2} \mu\left(\mathbf{r}^T(k)\mathbf{x}(k)\mathbf{x}^H(k)\right) +  \frac{1}{4} \mu\left(\mathbf{r}^T(k)\mathbf{x}(k)\mathbf{x}^T(k)\right)^* +\frac{1}{2} \mu n(k) \mathbf{x}^H(k) - \frac{1}{4} \mu \mathbf{x}^H(k)n^*(k)
\label{eq:convergence:CliveQLMS:v}
\end{equation}
Since QLMS is strictly linear, it is  second order optimal only for proper signals. For a second order circular $\mathbf{x}(k)$ we have
$
\mathbf{x}(k)\mathbf{x}^T(k)= -\frac{1}{2}\mathbf{x}(k)\mathbf{x}^H(k)
$
(see Appendix \ref{app:xxT_proof}) allowing us to  arrive at
\begin{equation*}
\mathbf{r}^T(k+1) = \mathbf{r}^T(k) -  \frac{1}{2} \mu\left(\mathbf{r}^T(k)\mathbf{x}(k)\mathbf{x}^H(k)\right) -  \frac{1}{8} \mu\left(\mathbf{r}^T(k)\mathbf{x}(k)\mathbf{x}^H(k)\right)^* +\frac{1}{2} \mu n(k) \mathbf{x}^H(k) - \frac{1}{4} \mu \mathbf{x}^H(k)n^*(k)
\end{equation*}

Similarly to the analysis of HR-QLMS, after splitting $\mathbf{r}(k)$ into its real and imaginary components (as shown in Appendix \ref{App:convergence}) and taking the statistical expectation, we have
\begin{equation}
\boldsymbol{\omega}^T(k+1)= \boldsymbol{\omega}^T(k)[\mathbf{I} - \frac{3}{4}\mu\mathbf{R}_c]
\end{equation}
where
\[
\mathbf{R}_c=
\left( \begin{array}{cccc}
\frac{5}{6}[\mathbf{R}_{x}]_r & {\imath}\frac{1}{2}[\mathbf{R}_{x}]_{\imath} & {\jmath}\frac{1}{2}[\mathbf{R}_{x}]_{\jmath} & {\kappa}\frac{1}{2}[\mathbf{R}_{x}]_{\kappa} \\
{\frac{5}{6}\imath}[\mathbf{R}_{x}]_{\imath} & \frac{1}{2}[\mathbf{R}_{x}]_r  & {\kappa}\frac{1}{2}[\mathbf{R}_{x}]_{\kappa} & {\jmath}\frac{1}{2}[\mathbf{R}_{x}]_{\jmath} \\
{\frac{5}{6}\jmath}[\mathbf{R}_{x}]_{\jmath} & {\kappa}\frac{1}{2}[\mathbf{R}_{x}]_{\kappa} & \frac{1}{2}[\mathbf{R}_{x}]_r  & {\imath}\frac{1}{2}[\mathbf{R}_{x}]_{\imath} \\
{\frac{5}{6}\kappa}[\mathbf{R}_{x}]_{\kappa} & {\jmath}\frac{1}{2}[\mathbf{R}_{x}]_{\jmath} & {\imath}\frac{1}{2}[\mathbf{R}_{x}]_{\imath} & \frac{1}{2}[\mathbf{R}_{x}]_r \\
\end{array} \right)
\]
Therefore, the condition for convergence in the mean becomes
\begin{equation}
|\mathbf{I} - \frac{3}{4}\mu\mathbf{R}_c|<1
\end{equation}
     and the bound on the stepsize
      \begin{equation}
 0 < \mu <\frac{8}{3\lambda_{max}(\mathbf{R}_c)}
\end{equation}

Repeating the analysis from Appendix \ref{App:convergence_Proof}, we can show that the eigenvalues of the matrices $\mathbf{R}_c$ and $\mathbf{R}_{xx}$ are related by
\begin{eqnarray}
\lambda_{max}(\mathbf{R}_c) \approx \frac{5}{6}\lambda_{max}(\mathbf{R}_{x}) \label{convergence:QLMS:max}\\
\lambda_{min}(\mathbf{R}_c) \approx \frac{1}{2}\lambda_{min}(\mathbf{R}_{x}) \label{convergence:QLMS:min}
\end{eqnarray}

\noindent \textbf{Remark\#17}: From the expressions for the maximum eigenvalue of the correlation matrices (\ref{convergence:HRQLMS:max}), (\ref{convergence:QLMS:max}), all the three algorithms have very similar stability properties, with the QLMS offering a slightly wider stability range.

\noindent \textbf{Remark\#18}: From  the analysis of the corresponding minimum eigenvalues of the correlation matrices (\ref{convergence:HRQLMS:min}), (\ref{convergence:QLMS:min}), the IQLMS is governed by the smallest eigenvalue spread \big($\frac{\lambda_{max}(\mathbf{R_{x}})}{\lambda_{min}(\mathbf{R_{x}})}$\big) and thus exhibits the fastest convergence rate, followed by QLMS and HR-QLMS.

We have illustrated that a structural equivalence exists among the three forms of QLMS algorithms. The analysis in the mean shows that in addition to requiring fewer computational operations, the IQLMS also offers a faster convergence rate. In view of the advantages offered by the IQLMS and the fact that IQLMS keeps the  generic form of LMS and CLMS and inherits their mathematical tractability, we provide the MSE analysis for the IQLMS as it is most likely to be used in practical applications.

\subsection{ Mean Square Error Performance of the Strictly Linear IQLMS for Noncircular Data}
   For generality, our analysis is performed for second order noncircular inputs, and naturally simplifies into the corresponding analysis for circular signals. The analysis evaluates the mean square error (MSE), defined as
   \vspace{0mm}
   \begin{equation}
  MSE = \lim_{k \rightarrow \infty} E[|e(k)|^2]
  \end{equation}
  for
  $
  e(k) = d(k) - \mathbf{w}^H(k)\mathbf{x}(k)
  $
  where $d(k)$ is an improper teaching signal. This is achieved using the approach proposed in \cite{Sayed_MSEperformance_conservationOfEnergy} \cite{Che_NonlinearAdaptive}, based on the energy conservation of filter weights at each iteration. To this end, we introduce the a priori error $e_a(k)$ and a  posteriori error  $e_p(k)$ as
  \vspace{0mm}
    \begin{equation*}
  e_a(k) = \mathbf{r}^H(k)\mathbf{x}(k) \hspace{20mm} e_p(k) = \mathbf{r}^H(k+1)\mathbf{x}(k)
  \vspace{0mm}
  \end{equation*}
  where $\mathbf{r}(k)= \mathbf{w}(k)- \mathbf{w}_o$ is the weight error vector. The error $e(k)$ can now be written as (see Appendix \ref{App:MSE_for_IQLMS})
  \vspace{-2mm}
  \begin{equation}
  e(k)=e_a(k) + \mathbf{w}^{cH}\mathbf{x}^a(k) + n(k)
  \label{MSE:IQLMS:e(k):ea(k)}
  \vspace{-2mm}
  \end{equation}
  and the MSE  as
  \vspace{-2mm}
  \begin{eqnarray}
MSE = EMSE + \lim_{k \rightarrow \infty}E\left[\| \mathbf{w}^{cH}\mathbf{x}^c(k)\|^2\right] + \sigma^2
 \label{eq:MSE:IQLMS:MSEexpression}
 \vspace{-2mm}
\end{eqnarray}
where $\sigma^2$ is the noise variance and EMSE is the excess mean square error, given by
\vspace{0mm}
 \begin{equation}
 EMSE = \lim_{k \rightarrow \infty} E[\|e_a(k)\|^2]
 \end{equation}
 while
 $
\mathbf{x}^c(k) = [\mathbf{x} ^{{\imath}T}(k),\mathbf{x} ^{{\jmath}T}(k) , \mathbf{x} ^{{\kappa}T}(k) ]^T$,
$\mathbf{w}^{c} = [ \mathbf{v}^T_o  ,\mathbf{g}^T_o  ,\mathbf{h}^T_o]^T$.
  The expression for the conservation of energy of the weights (see Appendix \ref{app:conservation_of_energy}) is given by
\begin{equation*}
\| \mathbf{w}(k+1)\|^2 + \frac{\|e_a(k)\|^2}{\|\mathbf{x}(k)\|^{2}}=
\| \mathbf{w}(k)\|^2 + \frac{\|e_p(k)\|^2} {\|\mathbf{x}(k)\|^{2}}
\end{equation*}
Upon taking the expectation of both sides, for the limit $k \rightarrow \infty$, we have
\begin{equation*}
E\left[\frac{\|e_a(k)\|^2}{\|\mathbf{x}(k)\|^{2}}\right]=
E\left[\frac{\|e_p(k)\|^2} {\|\mathbf{x}(k)\|^{2}}\right]
\end{equation*}
and upon substituting for $e_p(k) = e_a(k) - \mu e(k)\|\mathbf{x}(k)\|^2$ we arrive at
\vspace{0mm}
\begin{equation}
\mu^2E\left[\|\mathbf{x}(k)\|^2\|e(k)\|^2\right]= \mu E\left[e_a(k)e^*(k)\right] + \mu E \left[e(k) e_a^*(k)\right]
\label{eq:MSE:linear:e(k)}
\vspace{0mm}
\end{equation}
 Using the result in (\ref{eq:MSE:e(k) to e_a(k)}) and assumption A.1 in Appendix \ref{App:MSE_for_IQLMS}
 \vspace{0mm}
 \begin{equation}
 2\mu E[\|e_a(k)\|^2] = \mu^2 E\left[ \|\mathbf{x}(k)\|^2 \|e_a(k) \|^2\right]  + \mu^2Tr(\mathbf{R}_{x})  \sigma^2 + \mu^2 E \big[\|\mathbf{x}(k)\|^2  \| \mathbf{w}^{cH}\mathbf{x}(k) \|^2 \big]
 \label{eq:MSE:linear:final}
 \vspace{0mm}
 \end{equation}
 where $E\left[ \|\mathbf{x}\|^2 \right]=Tr(\mathbf{R}_{x})$. We can now obtain an expression for  $E[\|e_a(k)\|^2]$ under two different conditions:
\begin{itemize}
  \item For a small step size $\mu$ the term  $E\left[ \|\mathbf{x}(k)\|^2 \|e_a(k) \|^2\right]$ in  (\ref{eq:MSE:linear:final})  becomes negligible and  the expression for the EMSE can be written as
       \begin{equation}
  E[\|e_a(k)\|^2] = \frac{\mu}{2}Tr(\mathbf{R}_{x})  \sigma^2 + \mu^2 E \big[\|\mathbf{x}(k)\|^2  \| \mathbf{w}^{cH}\mathbf{x}(k) \|^2 \big]
    \label{eq:MSE:IQLMS:EMSE:small}
 \end{equation}
  \item For large values of $\mu$ where we cannot neglect the term $E\left[ \|\mathbf{x}(k)\|^2 \|e_a(k) \|^2\right]$, we make the usual independence assumption that $\|\mathbf{x}\|^2 $ is statistically independent to $\|e_a(k) \|^2$, giving
 \begin{equation}
  E[\|e_a(k)\|^2] = \frac{\mu Tr(\mathbf{R}_{x})  \sigma^2 + \mu^2 \big[ \|\mathbf{x}\|^2  \| \mathbf{w}^{cH}\mathbf{x}(k) \|^2 \big]}{2- \mu Tr(\mathbf{R}_x)}
  \label{eq:MSE:IQLMS:EMSE:big}
 \end{equation}
 \end{itemize}
 \subsection{ Mean Square Error Performance of the Widely Linear WL-IQLMS for Noncircular Data}
 The analysis for the MSE performance of the WL-IQLMS is similar to that carried out for the IQLMS, the difference arising from the relationship between $e(k)$ and $e_a(k)$, given by
 \vspace{0mm}
 \begin{equation}
 e(k)=e_a(k) + n(k)
 \vspace{0mm}
 \end{equation}
 that is, the term $\mathbf{w}^{cH}\mathbf{x}^c$ vanishes, and
\vspace{0mm}
 \begin{equation}
 MSE= EMSE + \sigma^2
 \label{eq:MSE:WLIQLMS:MSEexpression}
 \vspace{0mm}
 \end{equation}
\noindent \textbf{Remark\#19:}
 Comparing the MSE of WL-QLMS in (\ref{eq:MSE:WLIQLMS:MSEexpression}) to the MSE expression for the IQLMS in (\ref{eq:MSE:IQLMS:MSEexpression}), we observe that the term $\| \mathbf{w}^{cH}\mathbf{x}^c (k)\|^2 $ is always greater or equal to zero, and therefore subject on the EMSE for the IQLMS being smaller or equal than the EMSE for the WL-IQLMS, the MSE of the WL-IQLMS will always be smaller or equal than that of the  IQLMS.

The expressions for the EMSE of the WL-QLMS are thus given by:
 \begin{itemize}
  \item For a small step size
        \begin{equation}
  E[\|e_a(k)\|^2] = \frac{\mu}{2}Tr(\mathbf{R}_{x})  \sigma^2
  \label{eq:MSE:WLIQLMS:EMSE:small}
 \end{equation}
  \item For large values of the stepsize
 \begin{equation}
  E[\|e_a(k)\|^2] = \frac{\mu Tr(\mathbf{R}_{x})  \sigma^2 }{2- \mu Tr(\mathbf{R}_{x})}
  \label{eq:MSE:WLIQLMS:EMSE:big}
 \end{equation}
 \end{itemize}

 \noindent \textbf{Remark\#20:}
 Comparing (\ref{eq:MSE:WLIQLMS:EMSE:small}) and (\ref{eq:MSE:WLIQLMS:EMSE:big}) with  the EMSE expressions in (\ref{eq:MSE:IQLMS:EMSE:small})-(\ref{eq:MSE:IQLMS:EMSE:big}), we observe that for circular signals, where $\mathbf{w}^c=\mathbf{0}$, the EMSE of the WLIQLMS is four times greater than that of the IQLMS, though both tend to $0$ as the step size $\mu$ approaches 0. It follows from this that if the underlying model is know to be strictly, then the IQLMS should be used.
 From Remark \#19, it then  follows that the MSE of IQLMS is always smaller than that of WL-IQLMS if the signal is circular.

For noncircular signals, we observe that while the term containing the noise variance is four times larger for the WLIQLMS than the IQLMS, we must also consider the effect of the second term in the IQLMS EMSE containing the inner product $\| \mathbf{w}^{cH}\mathbf{x}^c(k) \|^2 $. For all noncircular signals, other than those that are marginally noncircular, the latter term will be the dominant factor and hence the EMSE of the IQLMS will be larger than that of the WLIQLMS for noncircular signals.  From Remark \#12, it then  follows that the MSE of IQLMS is always greater than that of WL-IQLMS if the signal is noncircular.

  \noindent \textbf{Remark\#21:}
For two-dimensional signals the analysis for the MSE of IQLMS and WL-IQLMS simplifies respectively into that of the complex LMS and the complex widely linear CLMS (CLMS and ACLMS).

\section{Simulations \label{sec:simulations}}
  In order to experimentally validate the analysis, we next evaluated the QLMS, HR-QLMS and IQLMS and their widely linear extensions on both circular and non-circular quaternion-valued data.

\begin{figure}[h]
\begin{minipage}[b]{0.5\linewidth}
\centering
\includegraphics[width=\linewidth,height=6cm]{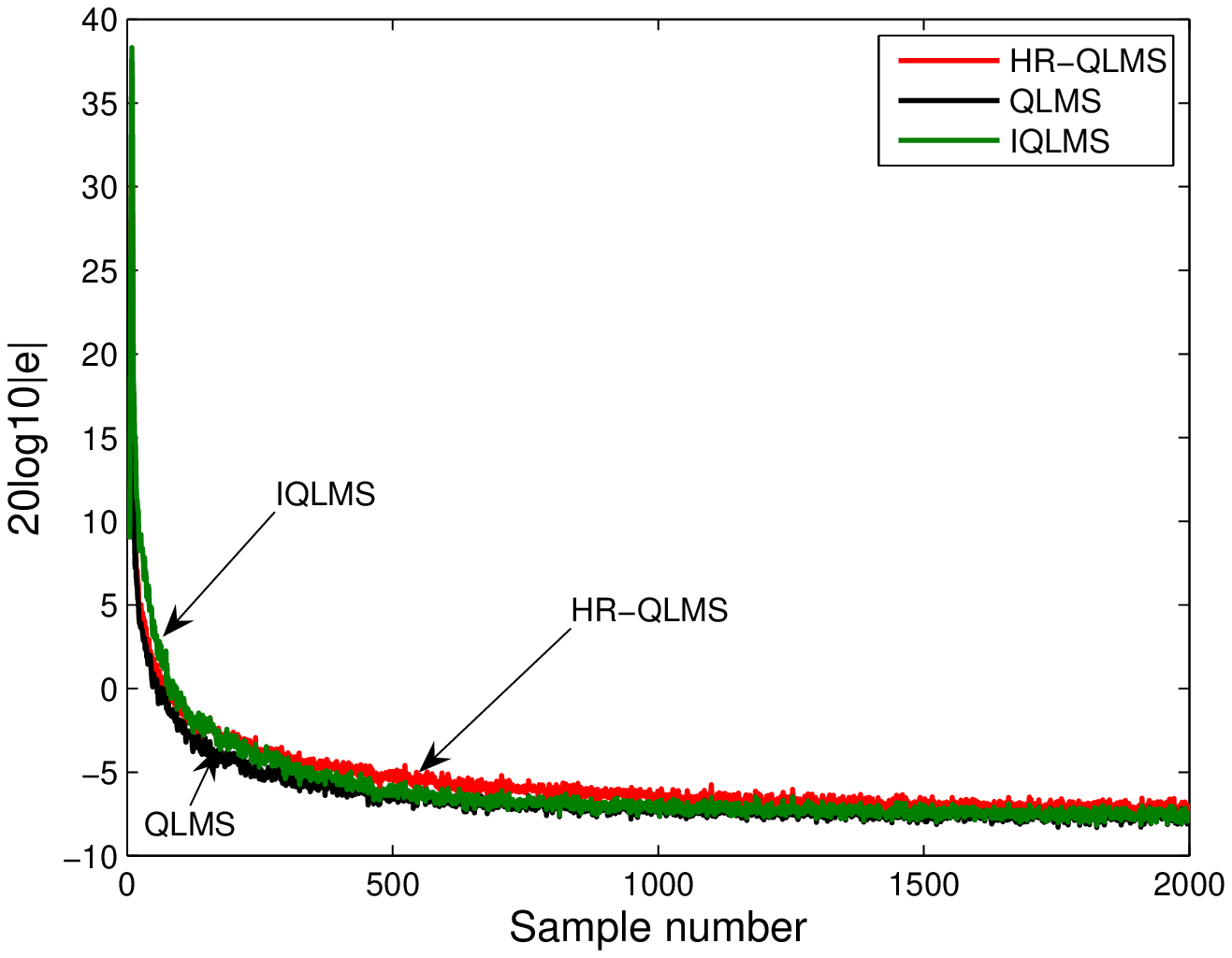}
\caption{Performance comparison of the QLMS, HR-QLMS and IQLMS on the  prediction of the circular AR(4) signal, for all the algorithms having the same step size.}
\label{Fig: AR4}
\end{minipage}
\hspace{0.5cm}
\begin{minipage}[b]{0.5\linewidth}
\centering
\includegraphics[width=\linewidth,height=6cm]{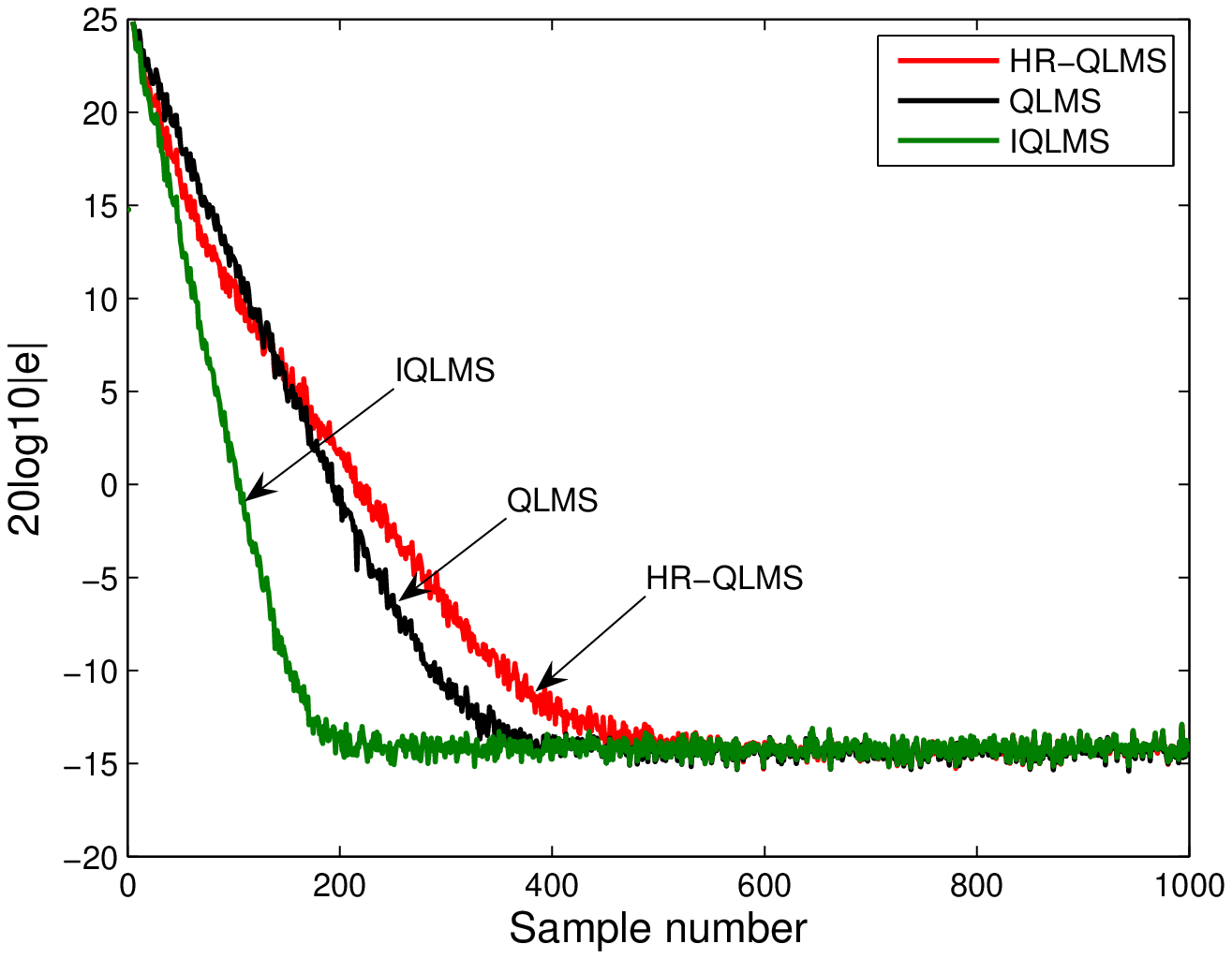}
\caption{Performance comparison of the QLMS, HR-QLMS and IQLMS on the prediction of the prediction of the circular MA(4) signal, for all the algorithms having the same step size.}
\label{Fig: MA4}
\end{minipage}
\end{figure}

 Fig. \ref{Fig: AR4} and Fig. \ref{Fig: MA4} show the performances of the HR-QLMS and IQLMS in the prediction setting on the circular AR(4) signal given by
 \vspace{0mm}
 \begin{equation}
 y(k)= 1.79y(k-1) - 1.85y(k-2) + 1.27y(k-3) - 0.41y(k-4) + n(k)
 \vspace{0mm}
 \end{equation}
 and the circular MA(4) signal
 \vspace{0mm}
 \begin{equation}
  y(k)= ax(k) + bx(k-1) + cx(k-2) + dx(k-3) + ex(k-4) + n(k)
  \vspace{0mm}
 \end{equation}
where $a,\dots,e$ are quaternion valued weights and $n(k)$ quadruply  white circular Gaussian  noise. For both benchmark models the driving noise $n(k)$ had a variance of 0.1. For the AR(4) signal, all the filters had a step size of $\mu = 0.08$ whereas for the MA(4) signal, all the filters had a step size of $\mu = 0.04$. The   experiments were repeated independently 100 times and  the performances were averaged to obtain the learning curves.      Observe from Fig. \ref{Fig: AR4} that for the circular AR(4) signal, IQLMS and QLMS offered marginally faster convergence than HR-QLMS, conforming with the analysis. In the steady state (for $k>10000$), all the three QLMS algorithms exhibited the same performance.  Fig. \ref{Fig: MA4} shows the performance of the MA(4) signal, illustrating the theoretical findings that all the three algorithms offer similar steady state performances, with the IQLMS exhibiting the fastest convergence. This also verifies the benefits of the I-gradient over the $\mathbb{HR}$ gradient.

\begin{figure}[h]
\centering
\includegraphics[width=12cm,height=10cm]{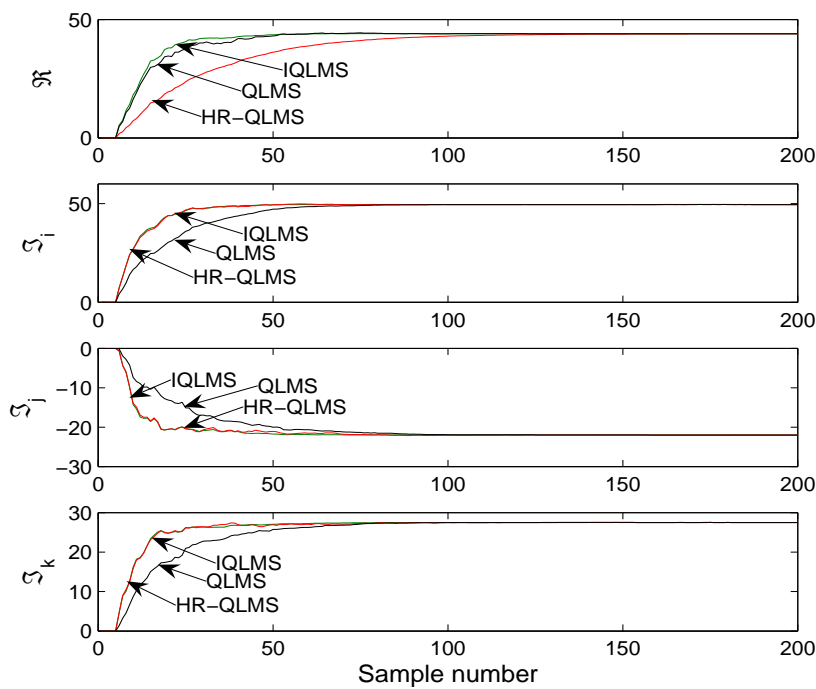}
\caption{Evolution of the adaptive quaternion valued weights of QLMS, HR-QLMS and IQLMS, for a single  trial prediction of the MA(4) signal.}
\label{Fig: weight_evolution}
\end{figure}

\begin{figure}[h]
\begin{minipage}[b]{0.5\linewidth}
\centering
\includegraphics[width=\linewidth,height=6cm]{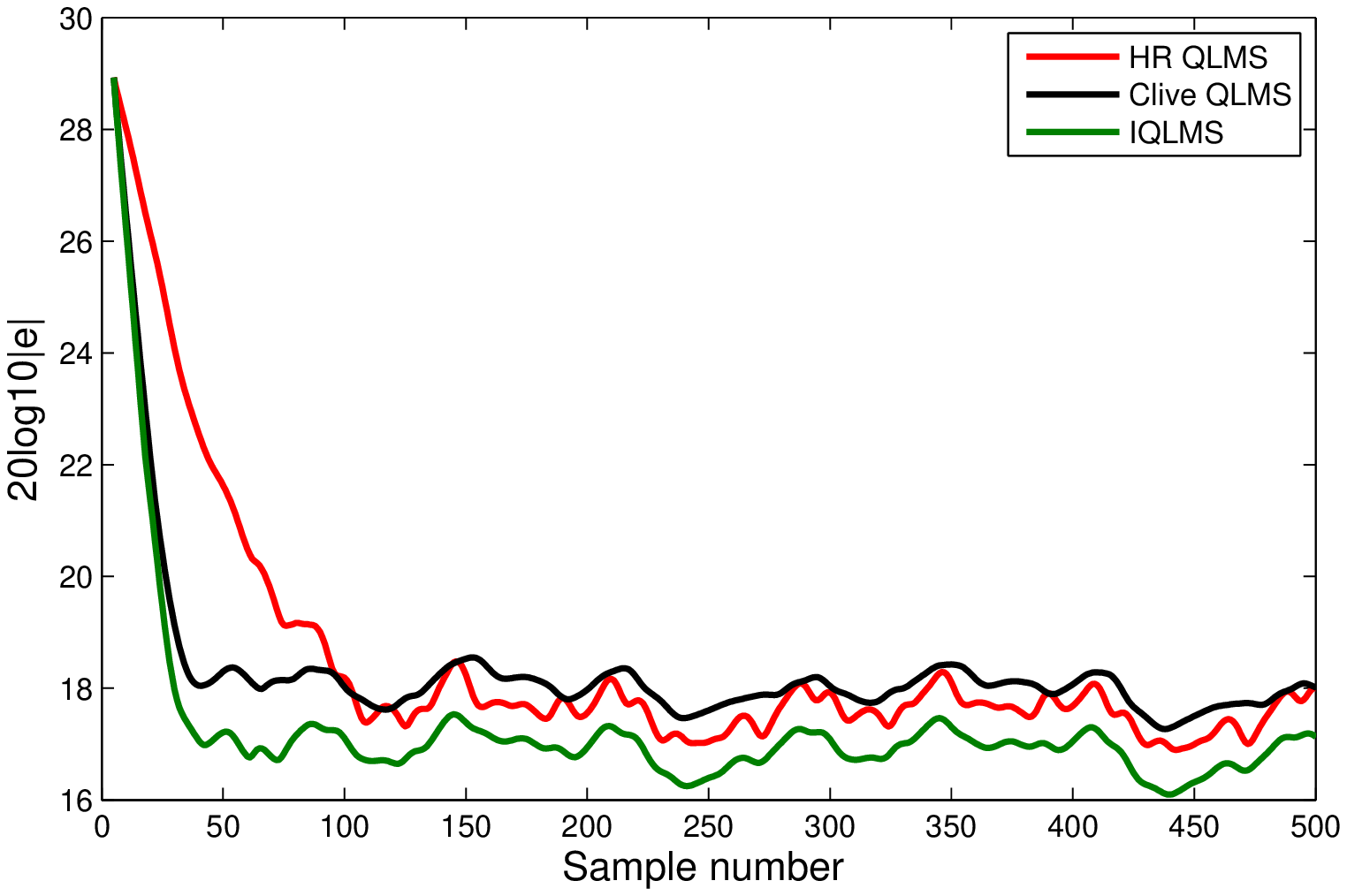}
\caption{Learning curves for the noncircular Lorenz signal, averaged over 100 experiments. }
\label{Fig: Lorenz3:QLMS}
\end{minipage}
\hspace{0.5cm}
\begin{minipage}[b]{0.5\linewidth}
\centering
\includegraphics[width=\linewidth,height=6cm]{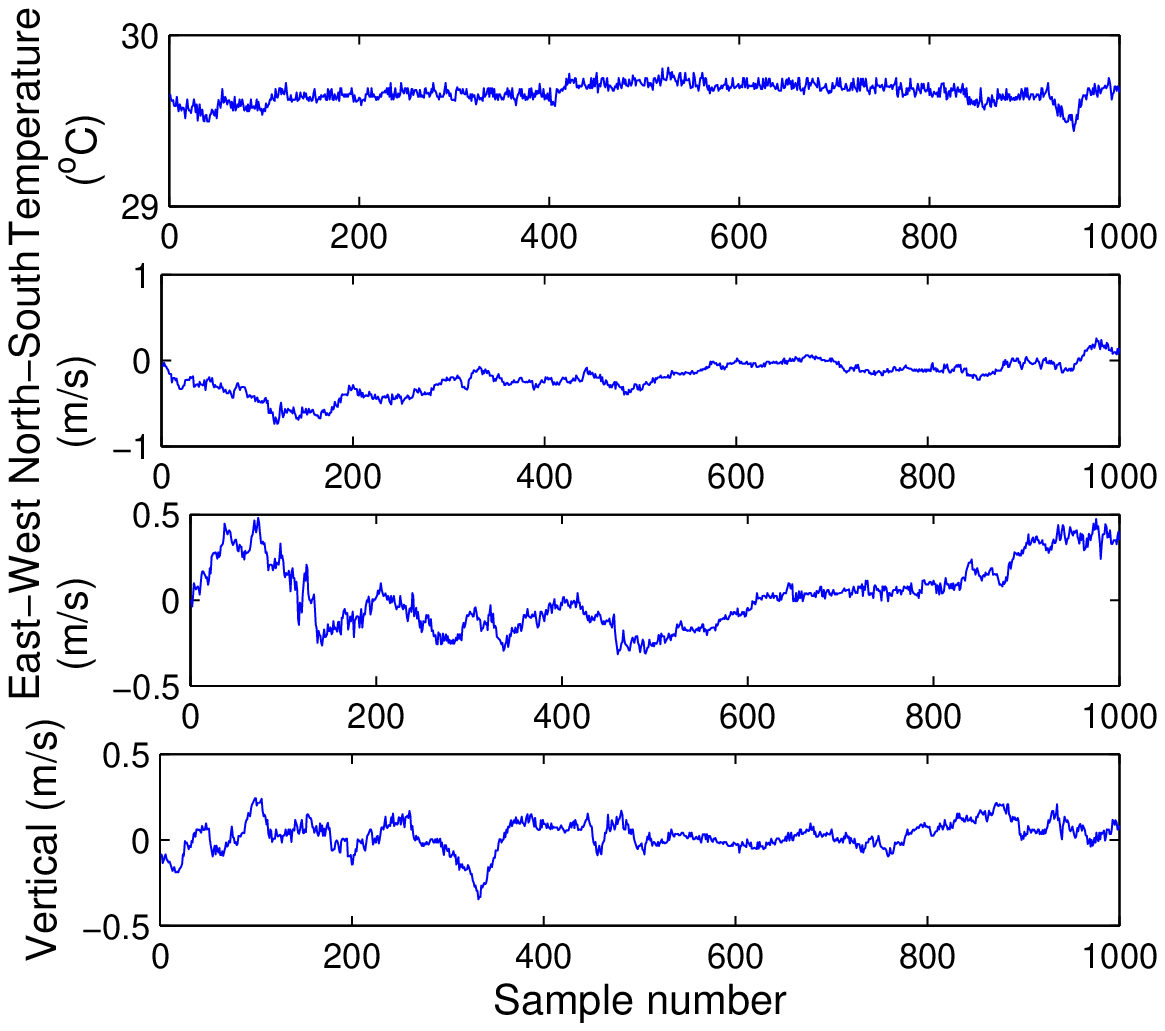}
\caption{Four-dimensional wind signal.  }
\label{Fig:wind_signal}
\end{minipage}
\end{figure}

 Fig. \ref{Fig: weight_evolution} shows the evolution of the real, $\imath-$, $\jmath-$, and $\kappa$-component of one of the filter weight vectors for the MA(4) signal. Verifying  (\ref{eq:equivalence3}) and (\ref{eq:equivalence4}), the IQLMS offered faster convergence for all the four data components. In addition, supporting the analysis, the HR-QLMS converged faster than the QLMS in all the three imaginary channels (with the QLMS converging faster in the real channel).
\\
Fig. \ref{Fig: Lorenz3:QLMS} compares the prediction performances of IQLMS, HR-QLMS  and QLMS  for the noncircular 3D Lorenz signal. A step size of $\mu = 2 \times 10^{-4}$ was used and the learning curves were averaged over 100 independent trials. Observe that the IQLMS achieved a slightly lower steady state error compared to QLMS and HR-QLMS, conforming with the analysis in Section \ref{Sec:Convergence}.

\begin{figure}[h]
\begin{minipage}[b]{0.5\linewidth}
\centering
\includegraphics[width=\linewidth,height=6cm]{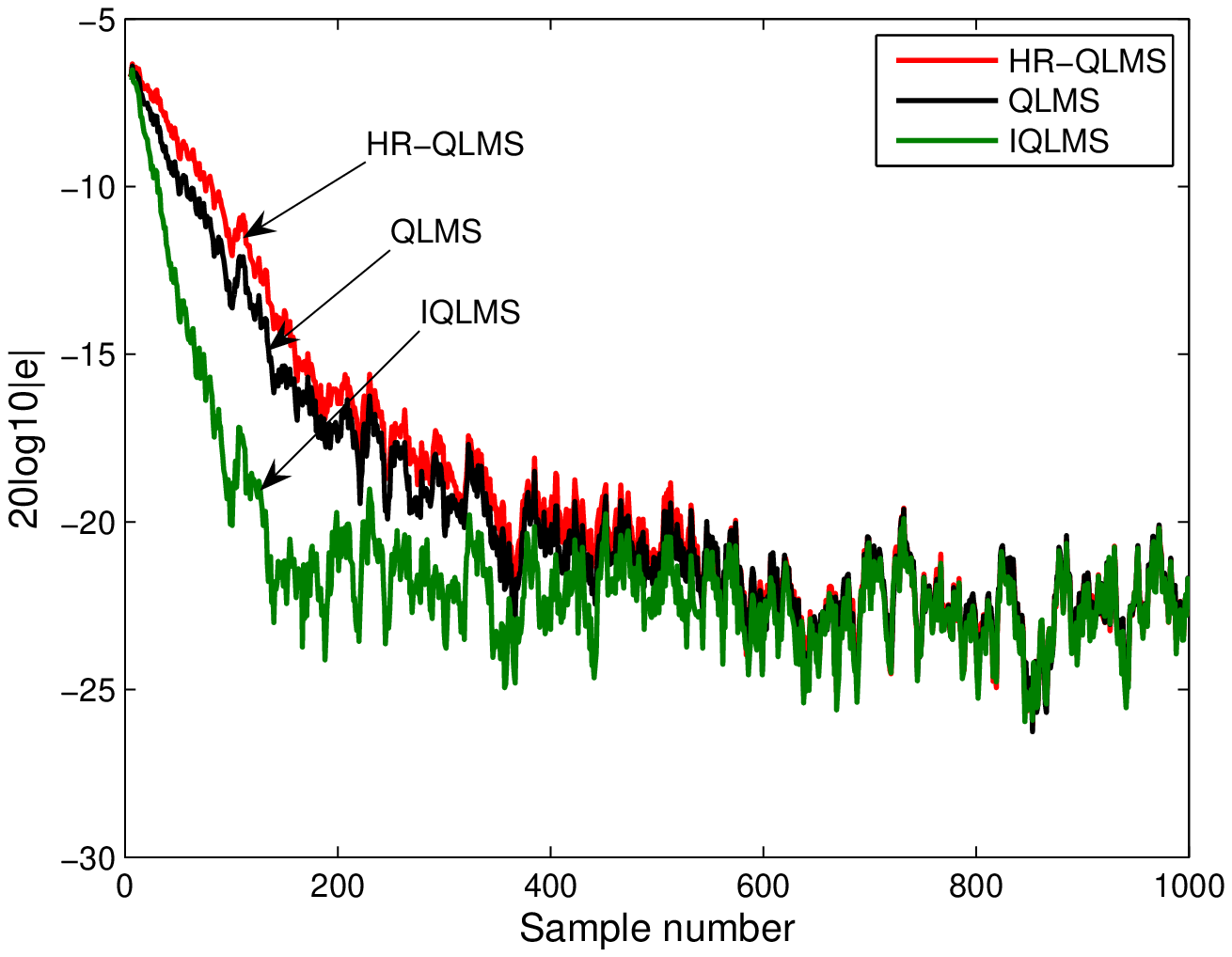}
\caption{Learning curves for the strictly linear algorithms on the prediction of the noncircular wind signal, averaged over 20 experiments. }
\label{Fig: windd_learning_curves}
\end{minipage}
\hspace{0.5cm}
\begin{minipage}[b]{0.5\linewidth}
\centering
\includegraphics[width=\linewidth,height=6cm]{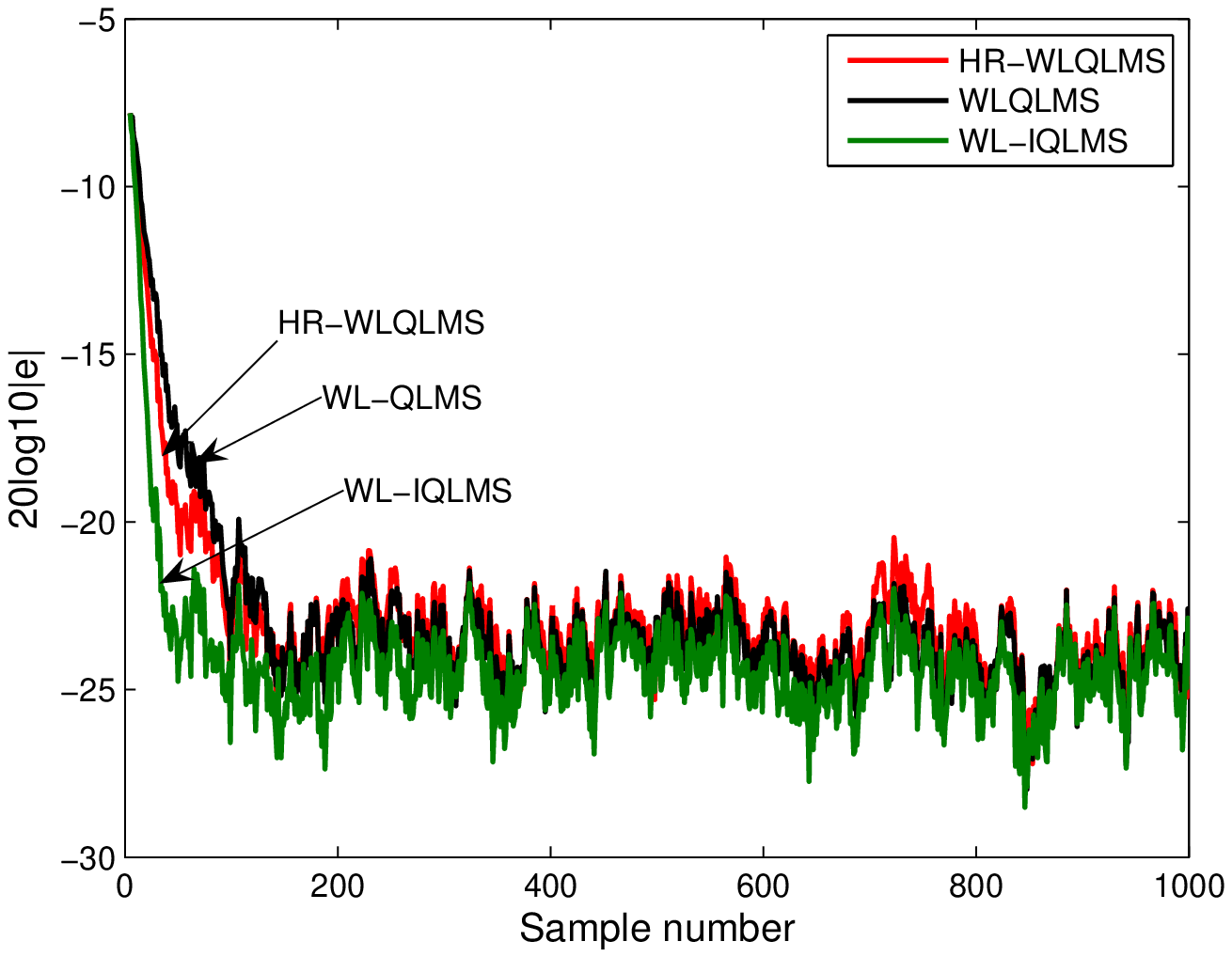}
\caption{Learning curves for the widely linear algorithms on the prediction of the  noncircular wind signal, averaged over 20 experiments. }
\label{Fig: windd_learning_curves_WL}
\end{minipage}
\end{figure}

Short term wind forecasting plays an important role in renewable energy, pollution modeling and aviation safety \cite{Manwell_WindEnergy_2002}. In this scenario, the QLMS, HR-QLMS and IQLMS were used to perform 10-step ahead prediction of a four dimensional wind signal\footnote{The wind data were recorded by Prof. K. Aihara and his team at the University of Tokyo, in an urban environment. The wind was initially sampled at 50 Hz, but resampled at 5 Hz for simulation purposes.}. The wind signal consists of three channels representing the wind speed measurement in the north-south, east-west and vertical direction and a fourth channel for the air temperature (see Fig. \ref{Fig:wind_signal}), and had a noncircularity degree\footnote{The noncircularity measure $r_s= \frac{E[xx^{\imath *}] + E[xx^{\jmath *}] + E[xx^{\kappa *}]}{3 E[xx^*]}$ where $r_s \in [0,1]$ and the value $r_s = 0$ indicates a circular source while $r_s=1$ indicates a highly noncircular source.} of  $r_s=0.49$. A step size of $2 \times 10^{-2}$ was used while the filters had a order of $N=4$. Fig. \ref{Fig: windd_learning_curves} shows the learning curves of the filters considered, averaged over 20 trials. Observe that all three filters achieved the same steady state performance and that the IQLMS converged considerably faster. The performances of the three widely linear filters are shown in Fig \ref{Fig: windd_learning_curves_WL}, using the same step size and filter order as for the strictly linear filters. Observe that all three filters achieved the same steady state performance while the WL-IQLMS exhibited the fastest convergence rate. Also note that, as expected from the noncircular nature of the wind signal, the widely linear filters achieved a better steady state performance than the corresponding strictly linear filters.
%

%
%

%

%

\section{Conclusions}
 This  work has resolved several ambiguities present in current stochastic gradient based quaternion valued adaptive filtering approaches, and has provided a rigorous unifying platform for the analysis and future developments in this field. We have shown that most inconsistencies in the treatment of quaternion valued algorithms arise due to  the noncommutativity of the quaternion product and the limitations of currently used quaternion gradients. To that end, we have introduced two new  gradient definitions in order to deal with the uncertainty related to the placement of imaginary units, and to enable direct differentiation of real functions of quaternion variables. The I-gradient has been shown to enable fastest convergence and to yield  a generic form of quaternion least mean square (QLMS) algorithms, both strictly linear and widely linear. We have illuminated that the class of IQLMS and WL-IQLMS algorithms are  natural extensions of LMS, CLMS and ACLMS, inheriting their mathematical compactness. For rigour, the mean and mean square convergence  analyses of the class of QLMS and widely linear QLMS (WL-QLMS) algorithms have been performed in a general widely linear setting and have been shown to naturally simplify into the corresponding analyses for the special cases (circular).  In this way, the rigorous treatment of quaternion gradients and noncommutativity of products presented in this work, together with convergence analysis and accounting for the noncircularity of distributions, has provided a  unifying platform for future development in quaternion valued adaptive filtering.

\appendices

 \section{The quaternion product rule  \label{app:product_rule}}
 The quaternion product rule is used to simplify the derivation of the quaternion valued algorithms proposed in this paper. To verify the validity of the product rule in the quaternion domain we rewrite from first principles the derivative of the product of two function
  \begin{eqnarray}
\frac{d(f(q)g(q))}{dq}&=& \lim_{dq\rightarrow 0} [f(q+dq)g(q+dq) - f(q)g(q)]dq^{-1} \nonumber \\
&=& \lim_{dq\rightarrow 0} [f(q+dq)g(q+dq) - f(q)g(q) + f(q)g(q+dq) - f(q)g(q+dq)]dq^{-1} \nonumber\\
&=& \lim_{dq\rightarrow 0}\Big[  f(q)[ g(q+dq) - g(q)] + [f(q+dq) - f(q) ]g(q+dq)  \Big] dq^{-1} \nonumber\\
&=& \lim_{dq\rightarrow 0}f(q)[g(q+dq) - g(q)]dq^{-1} + \lim_{dq\rightarrow 0}[f(dq+dq) - f(q) ]g(q+dq)dq^{-1} \nonumber\\
&=& f(q)\lim_{dq\rightarrow 0}[g(q+dq) - g(q)]dq^{-1} + \lim_{dq\rightarrow 0}[f(q+dq) - f(q) ]g(q)dq^{-1} \nonumber\\
&=&f(q)\frac{dg(q)}{dq} + \lim_{dq\rightarrow 0}[f(q+dq) - f(q) ]g(q)dq^{-1}
\label{two}
\end{eqnarray}
The second term does not seem to simplify into a quaternion derivative $\frac{\partial f(q)}{\partial q}$ due to the noncommutativity of the quaternion product. However, similar to the complex domain, we can choose the direction\footnote{The complex derivative  that form  the Cauchy Riemann equations can be equally expressed as
\begin{eqnarray*}
f'(z) &=& \frac{\partial u(x,y)}{\partial x} + \jmath  \frac{\partial v(x,y)}{\partial x} \\
f'(z)&=& \frac{\partial v(x,y)}{\partial y} - \jmath  \frac{\partial u(x,y)}{\partial y}
\end{eqnarray*} depending on whether we approach 0 from the $x$ or $y$ direction.} in which to approach $h \rightarrow 0$. Taking $dq$ along the path $hx$, where $h \rightarrow 0$ and $x$ is the real component of a quaternion, equation (\ref{two}) becomes
\begin{eqnarray*}
\frac{d(f(q)g(q))}{dq}&=&f(q)\frac{dg(q)}{dq} + \lim_{h\rightarrow 0}(f(q+dq) - f(q) )h^{-1}x g(q) \hspace{2 cm} \text{since $h^{-1}x$ is real}\\
\frac{d(f(q)g(q))}{dq}&=& f(q)\frac{dg(q)}{dq} + \frac{df(q)}{dq} g(q)
\end{eqnarray*}
which confirms the product rule and the validity of the analysis in this manuscript.

To illustrate  the validity of the above result on an example, we take the usual cost function (used throughout this paper) $f(q)=qq^*$. Using the product rule to simplify the calculation, the derivative becomes $\frac{\partial qq^*}{\partial q^*}=q - \frac{q^*}{2}= \frac{1}{2} q_r + \frac{3}{2}(q_r + \imath q_{\imath} + \jmath q_{\jmath} + \kappa q_{\kappa})$, while  multiplying out the product and taking the derivative we have $\frac{\partial qq^*}{\partial q^*}= -\frac{1}{2}q= \frac{1}{2}(q_r + \imath q_{\imath} + \jmath q_{\jmath} + \kappa q_{\kappa })$.

It then follows that for three-dimensional signals (where $q_r=0$) both derivatives point in the same direction (as exemplified by the sign of the $ q_{\imath}$ , $ q_{\jmath}$ and $q_{\kappa}$ being equal in both cases) and differ only in the amplitude of the gradient, $\frac{1}{2}$ vs $\frac{3}{2}$, which can be absorbed into the learning rate. For four-dimensional signals, all the  components also have the same sign, and will therefore reach the same minimum of the error surface. However, even after adjusting for the scaling factor, the magnitude of the real part will not scale equally and will therefore have a different convergence rate.  The I-gradient rectifies this issue and is shown to produce the same generic forms of stochastic gradient adaptive filters as in the other two division algebras, the reals and complex numbers. This makes it an enabling tool for future developments in this field, which has so far been upheld by the lack of rigorous and intuitive gradient.

\section{Duality between four-channel real valued LMS filter and the QLMS \label{app:duality}}
To show the equivalence between the WL-IQLMS and four-channel real LMS we compare the signal output and weight update of both filters and show that they are identical.
\subsection{The four-channel real valued LMS filter}

The output of the four-channel real valued LMS filter is given by:
\begin{equation}
\underbrace{\left[ \begin{array}{c}
y_r(k) \\
y_{\imath}(k) \\
y_{\jmath}(k)\\
y_{\kappa}(k) \\
\end{array} \right]}_{\mathbf{y}}= \frac{1}{4}
\underbrace{\left[ \begin{array}{cccc}
 \mathbf{w}^T_{11}(k) & \mathbf{w}^T_{12}(k) & \mathbf{w}^T_{13}(k) & \mathbf{w}^T_{14}(k) \\
 \mathbf{w}^T_{21}(k) &  \mathbf{w}^T_{22}(k) & \mathbf{w}^T_{23}(k) & \mathbf{w}^T_{24}(k)\\
 \mathbf{w}^T_{31}(k) &  \mathbf{w}^T_{32}(k) &  \mathbf{w}^T_{33}(k) & \mathbf{w}^T_{34}(k)\\
 \mathbf{w}^T_{41}(k) &  \mathbf{w}^T_{42}(k) &  \mathbf{w}^T_{43}(k) & \mathbf{w}^T_{44}(k) \end{array} \right]}_{\mathbf{W}}
\underbrace{\left[ \begin{array}{c}
\mathbf{x}_r(k) \\
\mathbf{x}_{\imath}(k) \\
\mathbf{x}_{\jmath}(k)\\
\mathbf{x}_{\kappa}(k) \\
\end{array} \right]}_{\mathbf{x}}
\label{dual1}
\end{equation}
where $\{\mathbf{w}_{11},\dots,\mathbf{w}_{44}\} \in \mathbb{R}^{N \times 1}$ are the 16 filter weights, $\{y_r,y_{\imath},y_{\jmath},y_{\kappa}\} \in \mathbb{R}^{1 \times 1}$ are the channel outputs  and $\{\mathbf{x}_r,\mathbf{x}_{\imath},\mathbf{x}_{\jmath},\mathbf{x}_{\kappa}\} \in \mathbb{R}^{N \times 1}$ are the input channels.
The error for each data channel is given by:
\begin{eqnarray}
e_{r}(k)&=& d_{r}(k)- y_r(k) \hspace{ 1cm}\nonumber
e_{{\imath}}(k)= d_{{\imath}}(k)- y_{\imath}(k) \nonumber \\
e_{{\jmath}}(k)&=& d_{{\jmath}}(k)- y_{\jmath}(k)  \hspace{ 1cm} \nonumber
e_{{\kappa}}(k)= d_{{\kappa}}(k)- y_{\kappa}(k)
\end{eqnarray}
and the objective function to be minimized
\begin{equation}
J=\frac{1}{2}(e^2_{r} + e^2_{q_{\imath}} + e^2_{q_{\jmath}} + e^2_{q_{\kappa}} )
\end{equation}
Based on (\ref{dual1}), the weight updates are obtained from $\mathbf{W}(k+1) = \mathbf{W}(k) - \mu  \nabla_\mathbf{W} J $, where the gradient $\nabla_\mathbf{W} J$ is given by
\begin{equation*}
\nabla_\mathbf{W} J = \left[ \begin{array}{cccc}
 \mathbf{e}_{r}(k) \mathbf{x}_r(k) & \mathbf{e}_{r}(k) \mathbf{x}_{\imath}(k)  & \mathbf{e}_{r}(k) \mathbf{x}_{\jmath}(k) & \mathbf{e}_{r}(k) \mathbf{x}_{\kappa}(k) \\
 \mathbf{e}_{\imath}(k) \mathbf{x}_r(k) & \mathbf{e}_{\imath}(k) \mathbf{x}_{\imath}(k)  & \mathbf{e}_{\imath}(k) \mathbf{x}_{\jmath}(k) & \mathbf{e}_{\imath}(k) \mathbf{x}_{\kappa}(k) \\
  \mathbf{e}_{\jmath}(k) \mathbf{x}_r(k) & \mathbf{e}_{\jmath}(k) \mathbf{x}_{\imath}(k)  & \mathbf{e}_{\jmath}(k) \mathbf{x}_{\jmath}(k) & \mathbf{e}_{\jmath}(k) \mathbf{x}_{\kappa}(k) \\
  \mathbf{e}_{\kappa}(k) \mathbf{x}_r(k) & \mathbf{e}_{\kappa}(k) \mathbf{x}_{\imath}(k)  & \mathbf{e}_{\kappa}(k) \mathbf{x}_{\jmath}(k) & \mathbf{e}_{\kappa}(k) \mathbf{x}_{\kappa}(k)  \end{array} \right]
\end{equation*}

\subsection{Duality between the four-channel real LMS and the WL-IQLMS}

The output of the WL-IQLMS is given by
\begin{equation}
y(k) = \mathbf{u}^T(k)\mathbf{q}(k) + \mathbf{v}^T(k)\mathbf{x}^i(k) + \mathbf{g}^T(k)\mathbf{x}^j(k) + \mathbf{h}^T(k)\mathbf{x}^k(k) = \mathbf{w}^T(k)\mathbf{x}^a(k)
\end{equation}
where $ \mathbf{x}^a(k) =    [\mathbf{x}^T(k) , \mathbf{x}^{iT}(k) , \mathbf{x}^{jT}(k) , \mathbf{x}^{kT}(k)]^T$ and $ \mathbf{w}(k) =    [\mathbf{u}^T(k) , \mathbf{v}^T(k), \mathbf{g}^T(k) , \mathbf{h}^T(k)]^T$

We can compare the output of the WL-IQLMS to that of the  four-channel real LMS filter by writing the output $y(k)$ of the WL-IQLMS in the same form as the output of the four-channel real LMS, that is, in terms of $y_r(k)$, $y_{\imath}(k)$, $y_{\jmath}(k)$ and $y_{\kappa}(k)$. In doing so we obtain
\begin{eqnarray*}
y_{r}(k) &=&( \underbrace{\mathbf{u}_{r}(k) + \mathbf{v}_{r}(k) + \mathbf{g}_{r}(k) + \mathbf{h}_{r}(k)}_{\mathbf{w}_{11}(k)} )^T \mathbf{x}_{r}
+ ( \underbrace{-\mathbf{u}_{\imath}(k) - \mathbf{v}_{\imath}(k) + \mathbf{g}_{\imath}(k) + \mathbf{h}_{\imath}(k)}_{\mathbf{w}_{12}(k)} )^T \mathbf{x}_{\imath} \\
&+&( \underbrace{-\mathbf{u}_{\jmath}(k) + \mathbf{v}_{\jmath}(k) - \mathbf{g}_{\jmath}(k) + \mathbf{h}_{\jmath}(k)}_{\mathbf{w}_{13}(k)} )^T \mathbf{x}_{\jmath}
+ ( \underbrace{-\mathbf{u}_{\kappa}(k) + \mathbf{v}_{\kappa}(k) + \mathbf{g}_{\kappa}(k) - \mathbf{h}_{\kappa}(k)}_{\mathbf{w}_{14}(k)} )^T \mathbf{x}_{\kappa} \\
\end{eqnarray*}
\begin{eqnarray*}
y_{\imath}(k) &=& ( \underbrace{\mathbf{u}_{r}(k) + \mathbf{v}_{r}(k) - \mathbf{g}_{r}(k) - \mathbf{h}_{r}(k)}_{\mathbf{w}_{21}(k)} )^T \mathbf{x}_{\imath}
+ ( \underbrace{\mathbf{u}_{\imath}(k) + \mathbf{v}_{\imath}(k) + \mathbf{g}_{\imath}(k) + \mathbf{h}_{\imath}(k)}_{\mathbf{w}_{22}(k)} )^T \mathbf{x}_{r} \\
&+& ( \underbrace{\mathbf{u}_{\jmath}(k) - \mathbf{v}_{\jmath}(k) - \mathbf{g}_{\jmath}(k) - \mathbf{h}_{\jmath}(k)}_{\mathbf{w}_{23}(k)} )^T \mathbf{x}_{\kappa}
+ ( \underbrace{ - \mathbf{u}_{\kappa}(k) + \mathbf{v}_{\kappa}(k) - \mathbf{g}_{\kappa}(k) + \mathbf{h}_{\kappa}(k)}_{\mathbf{w}_{24}(k)} )^T \mathbf{x}_{\jmath} \\
\end{eqnarray*}
\begin{eqnarray*}
y_{\jmath}(k) =& & ( \underbrace{\mathbf{u}_{r}(k) - \mathbf{v}_{r}(k) + \mathbf{g}_{r}(k) - \mathbf{h}_{r}(k)}_{\mathbf{w}_{31}(k)} )^T \mathbf{x}_{\jmath}
+ ( - \underbrace{\mathbf{u}_{\imath}(k) + \mathbf{v}_{\imath}(k) + \mathbf{g}_{\imath}(k) - \mathbf{h}_{\imath}(k)}_{\mathbf{w}_{32}(k)} )^T \mathbf{x}_{\kappa} \\
&+& ( \underbrace{\mathbf{u}_{\jmath}(k) + \mathbf{v}_{\jmath}(k) + \mathbf{g}_{\jmath}(k) + \mathbf{h}_{\jmath}(k)}_{\mathbf{w}_{33}(k)} )^T \mathbf{x}_{r}
+ ( \underbrace{\mathbf{u}_{\kappa}(k) + \mathbf{v}_{\kappa}(k) - \mathbf{g}_{\kappa}(k) - \mathbf{h}_{\kappa}(k)}_{\mathbf{w}_{34}(k)} )^T \mathbf{x}_{\imath} \\
\end{eqnarray*}
\begin{eqnarray*}
y_{\kappa}(k) =& & ( \underbrace{\mathbf{u}_{r}(k) - \mathbf{v}_{r}(k) - \mathbf{g}_{r}(k) + \mathbf{h}_{r}(k)}_{\mathbf{w}_{41}(k)} )^T \mathbf{x}_{\kappa}
+ ( \underbrace{\mathbf{u}_{\imath}(k) - \mathbf{v}_{\imath}(k) + \mathbf{g}_{\imath}(k) - \mathbf{h}_{\imath}(k)}_{\mathbf{w}_{42}(k)} )^T \mathbf{x}_{\jmath} \\
&+& ( \underbrace{-\mathbf{u}_{\jmath}(k) - \mathbf{v}_{\jmath}(k) + \mathbf{g}_{\jmath}(k) + \mathbf{h}_{\jmath}(k)}_{\mathbf{w}_{43}(k)} )^T \mathbf{x}_{\imath}
+ ( \underbrace{\mathbf{u}_{\kappa}(k) + \mathbf{v}_{\kappa}(k) + \mathbf{g}_{\kappa}(k) + \mathbf{h}_{\kappa}(k)}_{\mathbf{w}_{44}(k)} )^T \mathbf{x}_{r} \\
\end{eqnarray*}
By comparing the equations above with (\ref{dual1}), we can identify the filter coefficient vectors $\mathbf{w}_{11}(k)$ to $\mathbf{w}_{44}(k)$, and the two filters therefore have the same output.

To show that the weight update for the WL-IQLMS is identical to the weight update of the four-channel real LMS, consider the WL-IQLMS weight update, given by
\begin{equation}
\textbf{u}(k+1) = \textbf{u}(k) + \mu\left(\frac{1}{2}e(k)\mathbf{x}^*(k)   \right) \nonumber
\hspace{1cm}
\textbf{v}(k+1) = \textbf{v}(k) + \mu\left(\frac{1}{2}e(k)\mathbf{x}^{i^*}(k)   \right) \nonumber
\end{equation}
\begin{equation}
\textbf{g}(k+1) = \textbf{g}(k) + \mu\left(\frac{1}{2}e(k)\textbf{x}^{j^*}(k) \right) \nonumber
\hspace{1cm}
\textbf{h}(k+1) = \textbf{h}(k) + \mu\left(\frac{1}{2}e(k)\textbf{x}^{k^*}(k) \right)
\label{eq:duality:weight_update}
\end{equation}
We can  compare (\ref{eq:duality:weight_update}) with the weight update for the four-channel real LMS. For illustration, we consider  the real part of the weight updates, that is, $\mathbf{u}_{r}$, $\mathbf{v}_{r}$, $\mathbf{g}_{r}$ and $\mathbf{h}_{r}$ and show that the real part of the weight update $\Delta \mathbf{w}_{11}(k)$ is related by:
\begin{equation}
\Delta \mathbf{w}_{11}(k)= \frac{1}{2} \big( \Delta \mathbf{u}_{r}(k) + \Delta \mathbf{v}_{r}(k) + \Delta \mathbf{g}_{r}(k) + \Delta \mathbf{h}_{r}(k) \big)
\end{equation}
To show this, we can write the weight update in (9) in terms on the real and imaginary components, as
\begin{eqnarray*}
\mathbf{u}_{r}(k+1)=\mathbf{u}_{r}(k) + \mu \left(\frac{1}{2}e_{r}(k)\mathbf{x}_{r}(k) + \frac{1}{2}e_{\imath}(k)\mathbf{x}_{\imath}(k) + \frac{1}{2}e_{\jmath}(k)\mathbf{x}_{\jmath}(k) + \frac{1}{2}e_{\kappa}(k)\mathbf{x}_{\kappa}(k) \right)
\end{eqnarray*}
\begin{eqnarray*}
\mathbf{v}_{r}(k+1)=\mathbf{v}_{r}(k) + \mu \left(\frac{1}{2}e_{r}(k)\mathbf{x}_{r}(k) + \frac{1}{2}e_{\imath}(k)\mathbf{x}_{\imath}(k) - \frac{1}{2}e_{\jmath}(k)\mathbf{x}_{\jmath}(k) - \frac{1}{2}e_{\kappa}(k)\mathbf{x}_{\kappa}(k)\right)
\end{eqnarray*}
\begin{eqnarray*}
\mathbf{g}_{r}(k+1)=\mathbf{g}_{r}(k) + \mu \left(\frac{1}{2}e_{r}(k)\mathbf{x}_{r}(k) - \frac{1}{2}e_{\imath}(k)\mathbf{x}_{\imath}(k) + \frac{1}{2}e_{\jmath}(k)\mathbf{x}_{\jmath}(k) - \frac{1}{2}e_{\kappa}(k)\mathbf{x}_{\kappa}(k)\right)
\end{eqnarray*}
\begin{eqnarray*}
\mathbf{h}_{r}(k+1)=\mathbf{h}_{r}(k) + \mu \left(\frac{1}{2}e_{r}(k)\mathbf{x}_{r}(k) - \frac{1}{2}e_{\imath}(k)\mathbf{x}_{\imath}(k) - \frac{1}{2}e_{\jmath}(k)\mathbf{x}_{\jmath}(k) + \frac{1}{2}e_{\kappa}(k)\mathbf{x}_{\kappa}(k) \right)
\end{eqnarray*}
Observing that
\begin{equation}
\frac{1}{2} \left(\Delta \mathbf{u}_{r}(k) + \Delta \mathbf{v}_{r}(k) + \Delta \mathbf{g}_{r}(k) + \Delta \mathbf{h}_{r}(k)\right) = \frac{1}{2}\left(2e_{r}(k)\mathbf{x}_{r}(k)   \right)= \Delta \mathbf{w}_{11}(k)
\end{equation}
this illustrates that the four-channel LMS and WL-IQLMS are equivalent when the four-channel LMS has a step size twice as large as that of the WL-IQLMS.

\section{Proof for the equation (\ref{eq:convergence:HRQLMS:v}) \label{app:proof_convergence_HR_v}}
Substitute for $e(k)$ into (\ref{Eq:HRQLMS_Convergence:w(k+1)}) to obtain
\begin{eqnarray*}
\mathbf{w}(k+1) = \mathbf{w}(k) + \frac{1}{2}\mu \mathbf{w}^T_o \mathbf{x}(k) \mathbf{x}^*(k) + \frac{1}{2}\mu n(k)\mathbf{x}^*(k) - \frac{1}{2} \mu \mathbf{w}^T(k)\mathbf{x}(k)\mathbf{x}^*(k) \nonumber \\
- \frac{1}{4} \mu \mathbf{x}(k)\mathbf{x}^H(k)\mathbf{w}^*_o - \frac{1}{4} \mu \mathbf{x}(k) n^*(k) + \frac{1}{4} \mu \mathbf{x}(k) \mathbf{x}^H(k) \mathbf{w}^*(k)
\end{eqnarray*}
We can now subtract the optimal weight vector $\mathbf{w}_o$ from both sides  to obtain the weight error vector $\mathbf{r}(k)=\mathbf{w}(k)-\mathbf{w}_o(k)$ in the form
\begin{equation*}
\mathbf{r}(k+1) = \mathbf{r}(k) -  \frac{1}{2} \mu\left(\mathbf{r}^T(k)\mathbf{x}(k)\mathbf{x}^*(k)\right) +  \frac{1}{4} \mu\left( \mathbf{x}(k)\mathbf{x}^H(k) \mathbf{r}^*(k) \right) +\frac{1}{2} \mu n(k) \mathbf{x}^*(k) - \frac{1}{4} \mu \mathbf{x}(k)n^*(k)
\end{equation*}
The Hermitian transpose of the  term $\mathbf{x}(k)\mathbf{x}^H(k) \mathbf{r}^*(k)$ can be written as $ \left( \mathbf{x}(k)\mathbf{x}^H(k) \mathbf{r}^*(k) \right)^H =  \mathbf{r}^T(k)\mathbf{x}(k) \mathbf{x}^H(k)$ and noting also that $\mathbf{r}^T(k)\mathbf{x}(k)\mathbf{x}^*(k) = (\mathbf{r}^T(k)\mathbf{x}(k)\mathbf{x}^H(k))^T$, we have
\begin{equation*}
\mathbf{r}(k+1) = \mathbf{r}(k) -  \frac{1}{2} \mu\left(\mathbf{r}^T(k)\mathbf{x}(k)\mathbf{x}^H(k)\right)^T +  \frac{1}{4} \mu\left(\mathbf{r}^T(k)\mathbf{x}(k)\mathbf{x}^H(k)\right)^H +\frac{1}{2} \mu n(k) \mathbf{x}^*(k) - \frac{1}{4} \mu \mathbf{x}(k)n^*(k)
\end{equation*}
Taking the transpose of both sides gives
\begin{equation*}
\mathbf{r}^T(k+1) = \mathbf{r}^T(k) -  \frac{1}{2} \mu\left(\mathbf{r}^T(k)\mathbf{x}(k)\mathbf{x}^H(k)\right) +  \frac{1}{4} \mu\left(\mathbf{r}^T(k)\mathbf{x}(k)\mathbf{x}^H(k)\right)^* +\frac{1}{2} \mu n(k) \mathbf{x}^H(k) - \frac{1}{4} \mu \mathbf{x}^T(k)n^*(k)
\end{equation*}

\section{Expression for the convergence of the QLMS \label{App:convergence}}
Upon applying the statistical expectation operator to (\ref{eq:convergence:HRQLMS:v}) we obtain
\vspace{0mm}
\begin{equation}
E[\mathbf{r}^T(k+1)] = E[\mathbf{r}^T(k)] -  \frac{1}{2} \mu\left(E[\mathbf{r}^T(k)]\mathbf{R}_{x}\right) +  \frac{1}{4} \mu\left(E[\mathbf{r}^T(k)]\mathbf{R}_{x}\right)^*
\vspace{0mm}
\end{equation}
For $\mathbf{r}^T(k)=\mathbf{r}^T_r(k) + \imath\mathbf{r}^T_{\imath}(k) + \jmath\mathbf{r}^T_{\jmath}(k) + \kappa\mathbf{r}^T_{\kappa}(k)$, we can write
\begin{eqnarray*}
E[\mathbf{r}^T_r(k+1)] &=& E[\mathbf{r}^T_r(k)] - \frac{1}{2}\mu \left[E[\mathbf{r}^T_r(k)][\mathbf{R}_{x}]_r + {\imath}E[\mathbf{r}^T_{\imath}(k)]  {\imath}[\mathbf{R}_{x}]_{\imath}   + {\jmath}E[\mathbf{r}^T_{\jmath}(k)] {\jmath}[\mathbf{R}_{x}]_{\jmath} \right. \\
&+&  \left. {\kappa}E[\mathbf{r}^T_{\kappa}(k)]  {\kappa}[\mathbf{R}_{x}]_r  \right]
 + \frac{1}{4}\mu \left[E[\mathbf{r}^T_r(k)][\mathbf{R}_{x}]_r + {\imath}E[\mathbf{r}^T_{\imath}(k)] {\imath}[\mathbf{R}_{x}]_{\imath}  \right.\\
 &+& \left.
   {\jmath}E[\mathbf{r}^T_{\jmath}(k)] {\jmath}[\mathbf{R}_{x}]_{\jmath}   + {\kappa}E[\mathbf{r}^T_{\kappa}(k)] {\kappa}[\mathbf{R}_{x}]_r  \right] \\
 &=& E[\mathbf{r}^T_r(k)]  -\frac{1}{4}\mu \left[E[\mathbf{r}^T_r(k)][\mathbf{R}_{x}]_r + {\imath}E[\mathbf{r}^T_{\imath}(k) ] {\imath}[\mathbf{R}_{x}]_{\imath}   + {\jmath}E[\mathbf{r}^T_{\jmath}(k)] {\jmath}[\mathbf{R}_{x}]_{\jmath} \right.  \\
 &+& \left. {\kappa}E[\mathbf{r}^T_{\kappa}(k)] \kappa[\mathbf{R}_{x}]_r  \right]
\end{eqnarray*}
and similarly
\begin{eqnarray*}
{\imath} E[\mathbf{r}^T_{\imath}(k+1)]
 &=& {\imath} E[\mathbf{r}^T_{\imath}(k)] -\frac{3}{4}\mu \left[E[\mathbf{r}^T_r(k)] {\imath}[\mathbf{R}_{x}]_{\imath} + {\imath}E[\mathbf{r}^T_{\imath}(k)] [\mathbf{R}_{x}]_r   + {\jmath}E[\mathbf{r}^T_{\jmath}(k)] {\kappa}[\mathbf{R}_{x}]_{\kappa}   \right. \\
 &+& \left. {\kappa}E[\mathbf{r}^T_{\kappa}(k)]  {\jmath}[\mathbf{R}_{x}]_{\jmath}  \right]
\end{eqnarray*}
\begin{eqnarray*}
{\jmath} E[\mathbf{r}^T_{\jmath}(k+1)]
 &=& {\jmath} E[\mathbf{r}^T_{\jmath}(k)] -\frac{3}{4}\mu \left[ E[\mathbf{r}^T_r(k)] {\jmath}[\mathbf{R}_{x}]_{\jmath} + {\imath} E[\mathbf{r}^T_{\imath}(k)]  {\kappa}[\mathbf{R}_{x}]_{\kappa}  + {\jmath} E[\mathbf{r}^T_{\jmath}(k)] [\mathbf{R}_{x}]_r   \right. \\
 &+& \left. {\kappa} E[\mathbf{r}^T_{\kappa}(k)] {\imath}[\mathbf{R}_{x}]_{\imath}    \right]
\end{eqnarray*}
\begin{eqnarray*}
{\kappa} E[\mathbf{r}^T_{\kappa}(k+1)]
 &=& {\jmath} E[\mathbf{r}^T_{\jmath}(k)] -\frac{3}{4}\mu \left[ E[\mathbf{r}^T_r(k)] {\kappa}[\mathbf{R}_{x}]_{\kappa} + {\imath} E[\mathbf{r}^T_{\imath}(k)]  {\jmath}[\mathbf{R}_{x}]_{\jmath}  + {\jmath} E[\mathbf{r}^T_{\jmath}(k)] {\imath}[\mathbf{R}_{x}]_{\imath}  \right. \\
 &+& \left. {\kappa} E[\mathbf{r}^T_{\kappa}(k)] [\mathbf{R}_{x}]_r    \right]
\end{eqnarray*}

\section{Convergence comparison of the QLMS and IQLMS \label{App:convergence_Proof}}
To compare the eigenvalues of the covariance matrix $\mathbf{R}_b$ within the HR-QLMS to those of the matrix $\mathbf{R}_{xx}$ within the IQLMS, we  re-write the weight error in (\ref{eq:IQLMS_error_update}) in terms of its $r-$, $\imath -$, $\jmath -$ and $\kappa -$parts of the weight error vector as
\begin{equation}
E[\boldsymbol{\omega}(k+1)] = (\mathbf{I} - \frac{3}{4}\mu \mathbf{R}_a)E[\boldsymbol{\omega}(k)]
\end{equation}
where ${\boldsymbol{\omega}}(k) = [\mathbf{v}^T_r(k) \text{ , } \imath\mathbf{v}^T_\imath(k) \text{ , } \jmath\mathbf{v}^T_\jmath(k) \text{ , } \kappa\mathbf{v}^T_\kappa(k)]^T $ and
\[
\mathbf{R}_a=
\left( \begin{array}{cccc}
[\mathbf{R}_{x}]_r & \imath[\mathbf{R}_{x}]_{\imath} & \jmath[\mathbf{R}_{x}]_{\jmath} & \kappa[\mathbf{R}_{x}]_{\kappa} \\
\imath[\mathbf{R}_{x}]_{\imath} & [\mathbf{R}_{x}]_r  & \kappa[\mathbf{R}_{x}]_{\kappa} & \jmath[\mathbf{R}_{x}]_{\jmath} \\
\jmath[\mathbf{R}_{x}]_{\jmath} & \kappa[\mathbf{R}_{x}]_{\kappa} & [\mathbf{R}_{x}]_r  & \imath[\mathbf{R}_{x}]_{\imath} \\
\kappa[\mathbf{R}_{x}]_{\kappa} & \jmath[\mathbf{R}_{x}]_{\jmath} & \imath[\mathbf{R}_{x}]_{\imath} & [\mathbf{R}_{x}]_r \\
\end{array} \right)
\]
Notice that the matrix $\mathbf{R}_b$ is identical to the matrix $\mathbf{R}_a$, except for the first column of the block matrices being scaled by 3, allowing us to express $\mathbf{R_b}$ in terms of  $\mathbf{R_a}$ as
\vspace{0mm}
\begin{equation}
\mathbf{R}_b=\mathbf{R}_a\mathbf{D}
\vspace{0mm}
\end{equation}
where $\mathbf{D}$ is the diagonal matrix $diag(diag(\frac{1}{3},..,\frac{1}{3})_N,\mathbf{I}_N,\mathbf{I}_N,\mathbf{I}_N)) \in \mathbb{H}^{4 \times N} $ and $\mathbf{I}_N$ is the identity matrix of size $N \times N$. We now make the assumption that the cross-covariance of the $r-$, $i-$, $j-$ and $k-$ part of signal vector $\mathbf{x}$ is either symmetric or zero. This condition may at first appear over-restrictive, but it is easy to show that every circular signal satisfies this condition \cite{Clive_Statistics_2010} \cite{Via_ICA}, for which both QLMS and IQLMS are optimal. Furthermore, there exists a much wider class of signals, other than circular, that satisfy this condition, since circular signals must also satisfy the following conditions [17]:
 \begin{eqnarray*}
 \mathbf{R}_{r r } &=& \mathbf{R}_{\imath \imath } = \mathbf{R}_{\jmath \jmath } = \mathbf{R}_{\kappa \kappa }  \hspace{20mm}
 \mathbf{R}_{r \imath } =   \mathbf{R}_{\jmath \kappa } \\
\mathbf{R}_{r \jmath } &=& -  \mathbf{R}_{\imath \kappa } \hspace{40mm}
\mathbf{R}_{r \kappa } =   \mathbf{R}_{\imath \jmath }
 \end{eqnarray*}
Making this assumption,  we have  $[\mathbf{R}_{x}]_\imath=[\mathbf{R}_{x}]_\jmath=[\mathbf{R}_{x}]_\kappa=\mathbf{0}$, and matrices $\mathbf{R_a}$ and $\mathbf{R_b}$ and $\mathbf{R_c}$ become block diagonal. Using the following:
 \begin{itemize}
   \item  The eigenvalues of a block diagonal matrix are the eigenvalues of its blocks;
   \item  Scaling a matrix by a factor $\alpha$ scales the eigenvalues by the same factor;
 \end{itemize}

we can further write
\begin{eqnarray}
\lambda_{max}(\mathbf{R_a})&=&\lambda_{max}(\mathbf{R}_{x}) \hspace{1cm} \lambda_{min}(\mathbf{R_a})=\lambda_{max}(\mathbf{R}_{x}) \\
\lambda_{max}(\mathbf{R_b})&=&\lambda_{max}(\mathbf{R}_{x}) \hspace{1cm} \lambda_{min}(\mathbf{R_b})=\frac{1}{3}\lambda_{max}(\mathbf{R}_{x}) \\
\lambda_{max}(\mathbf{R_c})&=&\frac{5}{6}\lambda_{max}(\mathbf{R}_{x}) \hspace{1cm} \lambda_{min}(\mathbf{R_c})=\frac{1}{2}\lambda_{max}(\mathbf{R}_{x})
\end{eqnarray}

\section{Proof for the equation (\ref{eq:convergence:CliveQLMS:v}) \label{app:proof_convergence:Clive_v}}
Substituting for $e(k)$ into (\ref{eq:equivalence1}), we obtain
\begin{eqnarray*}
\mathbf{w}(k+1) = \mathbf{w}(k) + \frac{1}{2}\mu \mathbf{w}^T_o \mathbf{x}(k) \mathbf{x}^*(k) + \frac{1}{2}\mu n(k)\mathbf{x}^*(k) - \frac{1}{2} \mu \mathbf{w}^T(k)\mathbf{x}(k)\mathbf{x}^*(k) \nonumber \\
- \frac{1}{4} \mu \mathbf{x}^*(k)\mathbf{x}^H(k)\mathbf{w}^*_o - \frac{1}{4} \mu \mathbf{x}^*(k) n^*(k) + \frac{1}{4} \mu \mathbf{x}^*(k) \mathbf{x}^H(k) \mathbf{w}^*(k)
\end{eqnarray*}
We can now subtract $\mathbf{w}_o$ to obtain the weight error vector $\mathbf{r}(k)=\mathbf{w}(k)-\mathbf{w}_o(k)$, in the form
\begin{equation*}
\mathbf{r}(k+1) = \mathbf{r}(k) -  \frac{1}{2} \mu\left(\mathbf{r}^T(k)\mathbf{x}(k)\mathbf{x}^*(k)\right) +  \frac{1}{4} \mu\left( \mathbf{x}^*(k)\mathbf{x}^H(k) \mathbf{r}^*(k) \right) +\frac{1}{2} \mu n(k) \mathbf{x}^*(k) - \frac{1}{4} \mu \mathbf{x}^*(k)n^*(k)
\end{equation*}
The Hermitian of the  term $\mathbf{x}^*(k)\mathbf{x}^H(k) \mathbf{r}^*(k)$ can be written as $ \left( \mathbf{x}^*(k)\mathbf{x}^H(k) \mathbf{r}^*(k) \right)^H =  \mathbf{r}^T(k)\mathbf{x}(k) \mathbf{x}^T(k)$ and noting also that $\mathbf{r}^T(k)\mathbf{x}(k)\mathbf{x}^*(k) = (\mathbf{r}^T(k)\mathbf{x}(k)\mathbf{x}^H(k))^T$, we have
\begin{equation*}
\mathbf{r}(k+1) = \mathbf{r}(k) -  \frac{1}{2} \mu\left(\mathbf{r}^T(k)\mathbf{x}(k)\mathbf{x}^H(k)\right)^T +  \frac{1}{4} \mu\left(\mathbf{r}^T(k)\mathbf{x}(k)\mathbf{x}^T(k)\right)^H +\frac{1}{2} \mu n(k) \mathbf{x}^*(k) - \frac{1}{4} \mu \mathbf{x}^*(k)n^*(k)
\end{equation*}
Taking the transpose of both sides gives
\begin{equation*}
\mathbf{r}^T(k+1) = \mathbf{r}^T(k) -  \frac{1}{2} \mu\left(\mathbf{r}^T(k)\mathbf{x}(k)\mathbf{x}^H(k)\right) +  \frac{1}{4} \mu\left(\mathbf{r}^T(k)\mathbf{x}(k)\mathbf{x}^T(k)\right)^* +\frac{1}{2} \mu n(k) \mathbf{x}^H(k) - \frac{1}{4} \mu \mathbf{x}^H(k)n^*(k)
\end{equation*}

\section{Proof for $E[\mathbf{xx}^T]=-\frac{1}{2}E[\mathbf{xx}^H]$ \label{app:xxT_proof}}
Take the transpose of  $\mathbf{x}= \frac{1}{2}(\mathbf{x}^{i*}+\mathbf{x}^{j*}+\mathbf{x}^{k*} -\mathbf{x}^*)$ and premultiply by $\mathbf{x}$ to obtain
\begin{equation*}
\mathbf{x}\mathbf{x}^T= \frac{1}{2}\mathbf{x}(\mathbf{x}^{i*}+\mathbf{x}^{j*}+\mathbf{x}^{k*} -\mathbf{x}^*)^T
= \frac{1}{2}\mathbf{x}(\mathbf{x}^{iH}+\mathbf{x}^{jH}+\mathbf{x}^{kH} -\mathbf{x}^H)
\vspace{0mm}
\end{equation*}
Taking the expectation of both sides gives
\vspace{0mm}
\begin{equation*}
\mathbf{C}_x=\frac{1}{2}(\mathbf{P}_x + \mathbf{S}_x + \mathbf{T}_x - \mathbf{R}_x)
\vspace{0mm}
\end{equation*}
where $\mathbf{C}_x= E[\mathbf{x}\mathbf{x}^T]$.
For circular signals, $\mathbf{P}_x=\mathbf{S}_x=\mathbf{T}_x=\mathbf{0}$, and therefore
\vspace{0mm}
\begin{equation*}
\mathbf{C}_x=-\frac{1}{2}\mathbf{R}_x
\end{equation*}

\section{MSE expression for the IQLMS \label{App:MSE_for_IQLMS}}
For an improper teaching signal $d(k)=\mathbf{w}^H_o\mathbf{x}(k) + \mathbf{v}^H_o\mathbf{x}^{\imath}(k) + \mathbf{g}^H_o\mathbf{x}^{\jmath}(k) + \mathbf{h}^H_o\mathbf{x}^{\kappa}(k) + n(k)$ the error can be written as
\begin{eqnarray}
e(k) &=& \mathbf{d}(k) - \mathbf{w}^H(k)\mathbf{x}(k) \nonumber \\
e(k) &=& \mathbf{w}^H_o\mathbf{x}(k) + \mathbf{v}^H_o\mathbf{x}^{\imath}(k) + \mathbf{g}^H_o\mathbf{x}^{\jmath}(k) + \mathbf{h}^H_o\mathbf{x}^{\kappa}(k) + n(k) -  \mathbf{w}^H(k)\mathbf{x}(k) \nonumber \\
e(k) &=& \mathbf{r}^H\mathbf{x}(k) + \mathbf{v}^H_o\mathbf{x}^{\imath}(k) + \mathbf{g}^H_o\mathbf{x}^{\jmath}(k) + \mathbf{h}^H_o\mathbf{x}^{\kappa}(k) + n(k) \nonumber \nonumber \\
e(k) &=& e_a(k) + \mathbf{w}^{aH}\mathbf{x}^a(k)  + n(k)
\label{eq:MSE:e(k) to e_a(k)}
\end{eqnarray}
where
\vspace{0mm}
\begin{equation*}
\mathbf{x}^c(k) = [\mathbf{x} ^{{\imath}T}(k) ,\mathbf{x} ^{{\jmath}T}(k) , \mathbf{x} ^{{\kappa}T}(k) ]^T \hspace{1cm}
\mathbf{w}^c = [ \mathbf{v}^T_o  ,\mathbf{g}^T_o  ,\mathbf{h}^T_o]^T
\vspace{0mm}
\end{equation*}
The mean square error can now be written as
\vspace{0mm}
\begin{equation*}
MSE = \lim_{k\rightarrow\infty} E \left[\| e(k)\|\right] = \lim_{k\rightarrow\infty} E \left[\| e_a(k) + n(k) + \mathbf{w}^{cH}\mathbf{x}^c(k) \|^2 \right]
\vspace{0mm}
\end{equation*}
Making the usual independent assumption that

\textbf{A.1} The terms  $e_a(k)$, $n(k)$ and $\mathbf{x}(k)$ are statistically independent and $\lim_{k \rightarrow \infty }E[e_a(k)]= E[n(k)] = 0$

 and the MSE takes the form
\begin{eqnarray*}
MSE &=& \lim_{k \rightarrow \infty} E\left[\|e_a(k)\|^2\right] + \lim_{k \rightarrow \infty}E\left[\| \mathbf{w}^{cH}\mathbf{x}^c\|^2\right] + \lim_{k \rightarrow \infty}E\left[\| n(k)\|^2\right]
= EMSE + \sigma^2
\end{eqnarray*}
where EMSE is the excess mean square error and is defined as
 \begin{equation*}
 EMSE = \lim_{k \rightarrow \infty} E[\|e_a(k)\|^2] + \lim_{k \rightarrow \infty}E\left[\| \mathbf{w}^{cH}\mathbf{x}^c\|^2\right]
 \end{equation*}

 \section{Proof for the conservation of weight error energy \label{app:conservation_of_energy}}
 Following the approach in  \cite{HyunChooolShin_Affine}, substitute for $\mathbf{r}(k) = \mathbf{w}_o - \mathbf{w}(k)$ into the IQLMS weight update (where the filter takes the form $\mathbf{w}^H\mathbf{x}$) to give
 \vspace{0mm}
  \begin{equation*}
  \mathbf{r}(k+1)= \mathbf{r}(k) - \mu\mathbf{x}(k)e^*(k)
  \vspace{0mm}
  \end{equation*}
  Taking the Hermitian of both sides we have
  \vspace{0mm}
  \begin{equation}
  \mathbf{r}^H(k+1)= \mathbf{r}^H(k) - \mu e(k)\mathbf{x}^H(k) \label{sec:affine:w tilde}
  \vspace{0mm}
  \end{equation}
  and by post-multiplying both sides by $\mathbf{x}(k)$ we obtain
  \vspace{0mm}
  \begin{equation}
    \mathbf{r}^H(k+1)\mathbf{x}(k)= \mathbf{r}^H(k)\mathbf{x}(k) - \mu e(k)\|\mathbf{x}(k)\|^2 \label{sec:affine:weight update 2}
    \vspace{0mm}
  \end{equation}
  The a priori error $e_a(k)$ and a posteriori error $e_p(k)$ can now be defined
\vspace{0mm}
  \begin{equation*}
  e_a(k) = \mathbf{r}^H(k)\mathbf{x}(k) \hspace{20mm} e_p(k) = \mathbf{r}^H(k+1)\mathbf{x}(k)
  \vspace{0mm}
  \end{equation*}
  Substitute into (\ref{sec:affine:weight update 2}) to give
\vspace{0mm}
\begin{eqnarray*}
e_p(k) &=& e_a(k) - \mu e(k)\|\mathbf{x}(k)\|^{2} \\
\mu e(k) &=& (e_a(k) - e_p(k))\|\mathbf{x}(k)\|^{-2}
\vspace{0mm}
\end{eqnarray*}
allowing us to  rewrite (\ref{sec:affine:w tilde}) as
\vspace{0mm}
\begin{equation*}
 \mathbf{r}^H(k+1)= \mathbf{r}^H(k) - \big(e_a(k) - e_p(k)\big)\|\mathbf{x}(k)\|^{-2}\mathbf{x}^H(k)
 \vspace{0mm}
\end{equation*}
Upon rearranging the terms above
\vspace{0mm}
\begin{equation*}
\mathbf{r}^H(k+1) + e_a(k)\|\mathbf{x}(k)\|^{-2}\mathbf{x}^H(k) = \mathbf{r}^H(k) + e_p(k)\|\mathbf{x}(k)\|^{-2}\mathbf{x}^H(k)
\vspace{0mm}
\end{equation*}
Evaluating the energy ($\| \cdot \|^2$) of both sides, the energy conservation relationship can be written as \cite{HyunChooolShin_Affine}
\begin{equation*}
\| \mathbf{w}(k+1)\|^2 + \frac{\|e_a(k)\|^2}{\|\mathbf{x}(k)\|^{2}}=
\| \mathbf{w}(k)\|^2 + \frac{\|e_p(k)\|^2} {\|\mathbf{x}(k)\|^{2}}
\end{equation*}

%
%

\bibliographystyle{ieeetr}
\bibliography{sample}

%
%
%

\end{document}